\documentclass[sn-mathphys-num]{sn-jnl}% Math and Physical Sciences Numbered Reference Style 
\usepackage{graphicx}%
\usepackage{multirow}%
\usepackage{amsmath,amssymb,amsfonts}%
\usepackage{amsthm}%
\usepackage{dsfont} 
\usepackage{mathrsfs}%
\usepackage[title]{appendix}%
\usepackage{xcolor}%
\usepackage{textcomp}%
\usepackage{manyfoot}%
\usepackage{booktabs}%
\usepackage{algorithm}%
\usepackage{algorithmicx}%
\usepackage{algpseudocode}%
\usepackage{listings}%
\theoremstyle{thmstyleone}%
\newtheorem{theorem}{Theorem}%  meant for continuous numbers
\newtheorem{proposition}[theorem]{Proposition}% 
\newtheorem{lemma}{Lemma}
\newtheorem*{notation}{Notation}
\newtheorem{corollary}{Corollary}

\theoremstyle{thmstyletwo}%
\newtheorem{remark}{Remark}%

\theoremstyle{thmstylethree}%
\newtheorem{definition}{Definition}%

\def \d {\mathrm{d}} % Element d'integration
\begin{document}

\title[ODE-heat system coupled] %Use the shortened version of the full title
{Null controllability of an ODE-heat system with coupled boundary and internal terms}
%{Null controllability of a coupled ODE-heat system internally and at the boundary}


\author[1]{\fnm{Idriss} \sur{Boutaayamou}}\email{d.boutaayamou@uiz.ac.ma}
\equalcont{These authors contributed equally to this work.}
\author*[2]{\fnm{Fouad} \sur{Et-tahri}}\email{fouad.et-tahri@edu.uiz.ac.ma}
\equalcont{These authors contributed equally to this work.}

\author[3,4]{\fnm{Lahcen} \sur{Maniar}}\email{maniar@uca.ma}
\equalcont{These authors contributed equally to this work.\\
\begin{center}
	\textbf{\textcolor{blue}{Dedicated to the memory of Professor Hammadi Bouslous}}
	\end{center} }

\affil[1]{\orgdiv{Lab-SIV}, \orgname{Polydisciplinary Faculty-Ouarzazate, Ibnou Zohr University}, \orgaddress{\street{BP 638}, \city{Ouarzazate}, \postcode{45000}, \country{Morocco}}}

\affil*[2]{\orgdiv{Lab-SIV}, \orgname{Faculty of Sciences-Agadir, Ibnou Zohr University}, \orgaddress{\street{B.P. 8106}, \city{Agadir},
		\country{Morocco}}}

\affil[3]{\orgdiv{Cadi Ayyad University}, \orgname{Faculty of Sciences Semlalia}, \orgaddress{\street{B.P. 2390}, \city{Marrakesh},
\country{Morocco}}}

\affil[4]{\orgdiv{University Mohammed VI Polytechnic}, \orgname{Vanguard Center}, \orgaddress{\city{Benguerir},  \country{Morocco}}}


\abstract{This paper is devoted to the theoretical and numerical analysis of the null controllability of a coupled ODE-heat system internally and at the boundary with Neumann boundary control.
	%of an ODE-heat system with coupled boundary and internal terms under Neumann boundary control. 
	First, we establish the null controllability of the ODE-heat with distributed control using Carleman estimates. Then, we conclude by the strategy of space domain extension. Finally, we illustrate the analysis with some numerical experiments.}




\maketitle
\tableofcontents
\section{Introduction and main results.}
Coupled ODE-heat systems are powerful tools for modeling complex phenomena where diffusion processes and local dynamics interact. They have applications in many fields, from thermal engineering to environmental and biomedical sciences. For example, the heat equation can model the behavior of a thermal sensor, and the aim may be to study the stabilization or controllability of a dimensional system by introducing sensor dynamics, see  \cite[chapter 15]{krstic2009delay}.
\par 
In this paper, we consider a situation where an ODE-heat system with coupled boundary and internal terms is to be controlled by a Neumann boundary control, we are interested in the null controllability as well as the numerical aspect of the theoretical results found of the following system:
%	The system we are dealing with is as follows:
\begin{equation}  \label{s0}
	\left\{
	\begin{aligned}
		&y_{t}-y_{xx} + a(x,t)y+b(x,t)z(t)=0
		& & \text {in}\; Q_{\ell}, \\
		& z^{\prime}(t)+c(t)z(t)-\kappa y_{x}(0,t)
		=0 & & \text {in}\;(0,T), \\
		& y(0,t)=\mu z(t)  & & \text {in}\;(0,T), \\
		& y_{x}(\ell,t) =u(t) & & \text {in}\;(0,T), \\
		& y(\cdot,0)=y_{0}  & & \text {in}\;(0,\ell), \\
		&z(0)=z_{0},   \\	
	\end{aligned}
	\right.
\end{equation}
where $T>0$ is a finite time, $\ell>0$ a fixed spatial end, $Q_{\ell}:=(0,\ell)\times (0,T)$. Here $z(t)\in \mathbb{R}$ is the state of ODE and $y(x,t)\in\mathbb{R}$ is the state of a heat equation which are coupled both internally by the heat flux $-\kappa y_{x}(0,t)$ and by the potential $b(x,t)$ and at the boundary via the coefficient $\mu$. The function $u\in L^{2}(0,T)$ acts as a
boundary control and is used to drive the state $(y,z)$ to $(0,0)$ at time $T$ from the initial state $(y_{0},z_{0})$. In the sequel, we will make the following assumptions.
\begin{itemize}
	\item The potentials are bounded:
	\begin{eqnarray}
		a, b\in L^{\infty}(Q_{\ell})\quad\mbox{and}\quad c\in L^{\infty}(0,T). \label{H1}
	\end{eqnarray}
	\item $\mu, \kappa\in\mathbb{R}$ are real such that 
	\begin{eqnarray}
		\mu\kappa>0. \label{H2}
	\end{eqnarray}
\end{itemize}
%	In this paper, we will assume that, the potential terms $a, b\in L^{\infty}(Q)$, $c\in L^{\infty}(0,T)$ are bounded and $\mu, \kappa\in\mathbb{R}$ such that $\mu\kappa>0$. 
The term $\mu\kappa$ appears as the coefficient of the heat flux at the boundary point in the new system, by making the change of state 
\begin{eqnarray}
	(y,z)\mapsto (Y,Z):=(y,\mu z). \label{chang}
\end{eqnarray}
The case $\mu\kappa<0$ represents a dissipative interaction, whereas $\mu\kappa=0$ is the non-interactive case. However, when $\mu\kappa>0$ we are in the presence of a reactive interaction. We will analyze system \eqref{s0} in case of a reactive interaction without making change \eqref{chang}. Then We look at the controllability properties of system \eqref{s0}, which will be the main topic of our paper. Our main finding is as follows.
\begin{theorem} \label{Main 1}
	Assume that assumptions \eqref{H1} and \eqref{H2} are satisfied. Then, for any $T>0$, the system \eqref{s0} is null controllable in time $T$. More precisely, for any $(y_{0}, z_{0})\in H^{1}(0,\ell)\times\mathbb{R}$ with the compatibility condition $y_{0}(0)=\mu z_{0}$, there exists a control $u\in H^{1/4}(0,T)$ such that the associated state $(y,z)$ of \eqref{s0} satisfies:
	\begin{eqnarray}
		y(\cdot,T)=0 \;\text{in}\; (0,\ell) \quad \mbox{and} \quad z(T)=0.  \label{Null}
	\end{eqnarray}
\end{theorem}
\begin{remark}
	Note that if $\kappa=0$ (non-interactive case), then the ODE on $z$ is decoupled from the first equation of \eqref{s0}. Consequently, for any initial data $(y_{0},z_{0})\in L^{2}(0,\ell)\times\mathbb{R}$ such that $z_{0}\neq 0$, we have $z(T)\neq 0$ by Cauchy Lipschitz theorem. Roughly speaking, in the case $\kappa\neq 0$, $y$ can be driven to $0$ by the control $u$ and $z$ can be driven to $0$ by the coupling term $\kappa y_{x}(0,t)$. 
	%\textcolor{red}{$\mu=0$...............}\\
\end{remark}
\par 
The problem of null controllability of the heat equation with Dirichlet conditions in the one-dimensional case was first proved by the method of moments by Hector Fattorini and David Russell, see \cite{fattorini1971exact} and \cite{tenenbaum2007new}. After this method of  Hector Fattorini and David Russell, the problem is solved independently by Gilles Lebeau, Luc Robbiano (see \cite{lebeau1995controle} and \cite{le2012carleman}) and Andrei Fursikov, Oleg Imanuvilov (see \cite{fursikov1996controllability}) with Carleman-type estimates.
Lebeau and Robbiano's approach consists in proving elliptic Carleman estimates. A spectral inequality can be deduced, i.e. a high-frequency control result. A method commonly referred to today as the Lebeau-Robbiano method then allows us to go from this high-frequency control result to null controllability results. We refer to \cite{le2012carleman} for details and to \cite{beauchard2018null, miller2010direct} for generalizations. The one-dimensional case was recently shown again using a backstepping approach by Jean-Michel Coron and Hoai-Minh Nguyen, see \cite{coron2017null}.
\par 
Null controllability for parabolic equations with different boundary condition scenarios has been widely studied in recent years \cite{fernandez2006null, fernandez2006global, fursikov1996controllability, khoutaibi2020null, maniar2017null} and for coupled PDE-PDE systems with dynamic boundary conditions, which deserve to be studied for their controllability and will eventually be addressed in a forthcoming paper, we cite \cite{berinde2023qualitative} and the references therein. The classical method for considering null controllability is the approach of Fursikov and Imanuvilov, which consists in establishing parabolic Carleman estimates. 
An observability inequality can be deduced and the latter is equivalent to null controllability. However, for boundary control
of coupled parabolic equations, it is not an easy task for obtaining the Carleman estimates. In this context, the method of moments and the backstepping approach have been used to overcome these difficulties, as in \cite{fernandez2010boundary} for parabolic coupled equations and \cite{zeng2024null} for ODE-heat coupled equations. For more details on the controllability of linear coupled parabolic systems, 
see the survey report \cite{ammar2011recent} and the references therein.
\par 
For system \eqref{s0}, we have used Fursikov and Imanuvilov's approach to a distributed control problem, and via the spatial domain extension strategy, we will derive controllability results for \eqref{s0}. Hence, to prove the null controllability of \eqref{s0}, we reformulate the boundary null controllability of \eqref{s0} as null controllability with distributed controls by extending the domain $(0, \ell)$ into $(0, L)$ with control acts in a region of $(\ell, L)$:
\begin{equation}  \label{s1}
	\left\{
	\begin{aligned}
		&y_{t}-y_{xx} + a(x,t)y+b(x,t)z(t)=\mathds{1}_{\omega}v
		& & \text {in}\; Q_{L}, \\
		& z^{\prime}(t)+c(t)z(t)-\kappa y_{x}(0,t)
		=0 & & \text {in}\;(0,T), \\
		& y(0,t)=\mu z(t)  & & \text {in}\;(0,T), \\
		& y_{x}(L,t) =0 & & \text {in}\;(0,T), \\
		& y(\cdot,0)=y_{0}  & & \text {in}\;(0,L), \\
		&z(0)=z_{0},   \\	
	\end{aligned}
	\right.
\end{equation}
where $L>\ell$ is a fixed spatial end, $\omega\subset (\ell, L)$ is a nonempty open subset, $\mathds{1}_{\omega}$ is the characteristic function of $\omega$. First, we examine the null controllability of \eqref{s1}, then we obtain the following main result.
\begin{theorem} \label{Main 2}
	Assume that assumption \eqref{H2} is satisfied. Then, for any $L>0$, $T>0$, any $a,b\in L^{\infty}(Q_{L})$, $c\in L^{\infty}(0,T)$, and any $\omega\subset (0,L)$ nonempty open subset, the system \eqref{s1} is null controllable in time $T$. More precisely, for any $(y_{0}, z_{0})\in L^{2}(0,L)\times\mathbb{R}$, there exists a control $v\in L^{2}(\omega_{T})$, where $\omega_{T}:=\omega\times (0,T)$ such that the associated state $(y,z)$ of \eqref{s1} satisfies :
	%			 $y(\cdot,T)=0$ in $(0,L)$ and $z(T)=0$. Moreover
	\begin{eqnarray*}
		y(\cdot,T)=0 \;\text{in}\; (0,L) \quad \mbox{and} \quad z(T)=0. 
	\end{eqnarray*}
	Moreover
	\begin{eqnarray}
		\int_{\omega_{T}}|v|^{2}\d x\d t\leq C\left(\int_{0}^{L}|y_{0}|^{2}\d x+ |z_{0}|^{2}\right) \label{control-state}
	\end{eqnarray}
	for some $C>0$.
	%			where $C_{op}$ is defined in \eqref{observability inequality1}.
\end{theorem}
For that, we will adopt Hilbert  uniqueness method, we refer to \cite{lions1988controlabilite}. Consequently, we will show the following observability inequality which is equivalent to the null controllability of \eqref{s1}: $\exists C>0,\;\forall 
(\varphi_{T},\rho_{T})\in L^{2}(0,L)\times\mathbb{R}$
\begin{eqnarray*}
	\quad\int_{0}^{L}|\varphi(x,0)|^{2}\d x+ |\rho(0)|^{2}\leq C\int_{\omega_T}|\varphi|^{2}\d x\d t, \label{observability inequality}
\end{eqnarray*}
where $(\varphi,\rho)$ is the solution of the following dual homogeneous backward problem of \eqref{s1} with respect to the inner product defined below in \eqref{scalar product}
\begin{equation}\label{s2}
	\left\{
	\begin{aligned}
		&-\varphi_{t}-\varphi_{xx} +a(x,t)\varphi=0
		& & \text {in}\;Q_L, \\
		& -\rho^{\prime}(t)+c(t)\rho(t)-\kappa\varphi_{x}(0,t) +\mu^{-1}\kappa\int_{0}^{L}b(x,t)\varphi(x,t)dx
		=0 & & \text {in}\;(0,T), \\
		& \varphi(0,t)=\mu \rho(t)  & & \text {in}\;(0,T), \\
		& \varphi_{x}(L,t) =0 & & \text {in}\;(0,T), \\
		& \varphi(\cdot,T)=\varphi_{T}  & & \text {in}\;(0,L), \\
		& \rho(T)=\rho_{T}.  	
	\end{aligned}
	\right.
\end{equation}
This form of the system \eqref{s2} is explained in Remark \ref{Remark 1}. The classic method for obtaining the observability inequality is to use a Carleman estimate. To our knowledge, a Carleman estimate for such a system has not been carried out in the literature. As a result, a new Carleman estimate has been proved,  see Proposition \ref{P4}.  The difficulty in this case is due to the nonlocality in space, and the boundary and internal couplings. Second, to obtain the controllability results of \eqref{s0}, we apply the spatial domain extension strategy as explained above.

\textbf{Structure of this paper.} This paper is organized as follows. In Section \ref{sec2}, we introduce the functional framework and the well-posedness of System \eqref{s1}. In Section \ref{sec3} we prove our Carleman estimate (Proposition \ref{P4}) and the main results of this paper (Theorems \ref{Main 1} and \ref{Main 2}). Section \ref{sec4} is devoted to numerical illustrations of the theoretical results found. Appendix \ref{Appendix A} is devoted to the proof of a few technical lemmas while Appendix \ref{Appendix B} is devoted to the proof of a duality relation.

\section{Semigroup generation and well-posedness} \label{sec2}
\subsection{General setting}
Let us first introduce some basic notation. Let $s>0$ be strictly positive real and $I$ an interval of $\mathbb{R}$, $L^2(I), L^{\infty}(I)$ and $H^{s}(I)$ are the usual spaces of Lebesgue and Sobolev for functions mapping from $I$ to $\mathbb{R}$. We write $\|\cdot\|_{L^{2}(I)}, \|\cdot\|_{L^{\infty}(I)}$ and $\|\cdot\|_{H^{s}(I)}$ to denote the standard norms on these spaces. We also denote by $D(I)$ the space of test functions on the interval $I$. For any Banach space $X$, the Bochner spaces of the functions mapping from $I$ to $X$ are denoted by $L^{2}(I;X), L^{\infty}(I;X)$ and $H^{s}(I;X)$, and the space $C(I;X)$ denotes the set of continuous functions mapping from $I$ to $X$. The spaces $L^{2}(0,T,L^{2}(I))$ and $L^{\infty}(0,T,L^{\infty}(I))$ can be identified, respectively, by $L^{2}(I\times (0,T))$ and $L^{\infty}(I\times (0,T))$.
The natural state space for our analysis is $$\mathcal{H}(I):=L^{2}\left(I; \d x\right)\times\mathbb{R},$$
where 
$\d x$ denotes the Lebesgue measure on $I$. This is a Hilbert space equipped with the following inner product
\begin{eqnarray}
	\left\langle (y,\alpha), (
	\varphi,\beta) \right\rangle_{\mathcal{H}(I)}:=\int_{I}y(x)\varphi(x)\d x+\frac{\mu}{\kappa}\alpha\beta,\quad (y,\alpha), (
	\varphi,\beta)\in \mathcal{H}(I). \label{scalar product}
\end{eqnarray}
This choice of inner product \eqref{scalar product} on $\mathcal{H}(I)$ has enabled us to symmetrize the main operator of our system. The Sobolev-type spaces compatible with our situation are defined by 
\begin{eqnarray*}
	\mathcal{H}^{s}(I) &:=& \left\{
	(y,\alpha)\in H^{s}(I)\times\mathbb{R}\colon y(0)=\mu \alpha \right\},\quad s\geq 1,
\end{eqnarray*}
equipped with the following inner product, is a Hilbert space
\begin{eqnarray*}
	\left\langle (y,\alpha), (\varphi,\beta) \right\rangle_{\mathcal{H}^{s}(I)} &:=& \langle y, \varphi \rangle_{H^{s}(I)}+\frac{\mu}{\kappa}\alpha\beta,\quad (y,\alpha),\; (\varphi,\beta)\in \mathcal{H}^{s}(I)\\
	&=& \sum_{j=0}^{s}\int_{I}y^{(j)}(x)\varphi^{(j)}(x)dx + \frac{\mu}{\kappa}\alpha\beta,
\end{eqnarray*}
where $s\in \mathbb{N}^{*}$, and for all $(y,\alpha)\in \mathcal{H}^{s}(I)$, we denote by $y^{\prime},\cdots, y^{(s)}$ the derivatives of $y$ in the distribution sense. For any $T>0$ strictly positive real, we also introduce the following energy space
\begin{eqnarray*}
	\mathcal{E}(I):= H^{1}(0,T;\mathcal{H}(I))\cap L^{2}(0,T;\mathcal{H}^{2}(I)),
\end{eqnarray*}
equipped with the following inner product, is a Hilbert space
\begin{equation*}
	\left\langle 	(y, \alpha), (
	\varphi, \beta) \right\rangle_{\mathcal{E}(I)} 
	:=\int_{I\times (0,T)}\left[y\varphi+ y_{t}\varphi_{t}+ y_{x}\varphi_{x}+ y_{xx}\varphi_{xx}\right]\d x\d t +\frac{\mu}{\kappa}\int_{0}^{T}\left[\alpha\beta+\alpha^{\prime}\beta^{\prime}\right]\d t,
\end{equation*}
where, for all $(y,\alpha)\in \mathcal{E}(I)$, we denote by $y_{t}, y_{x}$ and $y_{xx}$ the first and second derivatives of $y$ and $\alpha^{\prime}$ the first derivative with respect to $t$ of $\alpha$ in the sense of  distributions.
We denote by $(H^{1}(I))^{\prime}$ and $(\mathcal{H}^{1}(I))^{\prime}$ the dual spaces of $H^{1}(I)$ and $\mathcal{H}^{1}(I)$, respectively, and
\begin{eqnarray*}
	\mathcal{W}(I):=\left\{(y,\alpha)\in L^{2}(0,T;\mathcal{H}^{1}(I)) \colon (y_{t},\alpha^{\prime})\in L^{2}\left(0,T; (\mathcal{H}^{1}(I))^{\prime}\right)\right\}.
\end{eqnarray*}
We recall the usual continuous embedding
\begin{eqnarray*}
	\mathcal{W}(I)\hookrightarrow C([0,T];\mathcal{H}(I)).
\end{eqnarray*}
\begin{notation}
	In case $I=(0,L)$ for some $L>0$, unless otherwise specified, $Q_{L}$, $\mathcal{H}(0,L), \mathcal{H}^{s}(0,L), \mathcal{E}(0,L)$ and $\mathcal{W}(0,L)$ will be simply denoted by $Q, \mathcal{H}, \mathcal{H}^{s}, \mathcal{E}$ and $\mathcal{W}$ respectively. 
\end{notation}
\subsection{Well-posedness and regularity of the solution} 
In this section, we study the well-posedness and regularity of the following inhomogenuous system 
\begin{equation}\label{s3}
	\left\{
	\begin{aligned}
		&y_{t}-y_{xx} +a(x,t)y+ b(x,t)z(t)=f
		& & \text {in}\;Q, \\
		& z^{\prime}(t)+ c(t)z(t)-\kappa y_{x}(0,t) 
		=g & & \text {in}\;(0,T), \\
		& y(0,t)=\mu z(t)  & & \text {in}\;(0,T), \\
		& y_{x}(L,t) =0 & & \text {in}\;(0,T), \\
		& y(\cdot,0)=y_{0}  & & \text {in}\;(0,L), \\
		&z(0)=z_{0},   \\	
	\end{aligned}
	\right.
\end{equation} 
where $a,b\in L^{\infty}(Q), c\in L^{\infty}(0,T), \mu, \kappa\in\mathbb{R}$ such that $\mu\kappa>0$, $L>0$ and $(f,g)\in L^{2}(0,T;\mathcal{H})$.
The system \eqref{s3} can be written as an abstract Cauchy problem as follows
\begin{equation}  \label{acp}
	\left\{
	\begin{aligned}
		&\mathbf{U}^{\prime}(t)=(\mathcal{A}+\mathcal{B}(t))\mathbf{U}(t)+ \mathbf{F}(t)
		& & 0<t<T, \\
		&\mathbf{U}(0)=\mathbf{U}_{0},   
	\end{aligned}
	\right.
\end{equation}
where $\mathbf{U}(t):=(y(\cdot,t), z(t))$, $\mathbf{U}_0=(y_{0},z_{0})$, $\mathbf{F}(t)=(f(\cdot,t),  g(t))$ and the linear operators 
$$\mathcal{A} \colon D(\mathcal{A}) \subset \mathcal{H} \longrightarrow \mathcal{H}\quad\mbox{and}\quad  \mathcal{B}(t) \colon \mathcal{H} \longrightarrow \mathcal{H}$$
given by
\begin{eqnarray*}
	\mathcal{A}=
	\begin{pmatrix} \frac{d^{2}}{dx^{2}}  & 0\\ \kappa\frac{d}{dx}|_{x=0}  & 0 \end{pmatrix} \quad\mbox{and}\quad\mathcal{B}(t)=\begin{pmatrix} -a(\cdot,t)I_{L^{2}(0,L)}  & -b(\cdot,t)I_{\mathbb{R}}\\ 0  & -c(t)I_{\mathbb{R}} \end{pmatrix},
\end{eqnarray*}
where $I_{L^{2}(0,L)}$ and $I_{\mathbb{R}}$ are the identity operators, and the domains 
\begin{eqnarray*}
	D(\mathcal{A})=\left\{
	(y, \alpha)\in \mathcal{H}^{2}\quad : \quad y^{\prime}(L)=0\right\} \quad\mbox{and}\quad D(\mathcal{B}(t))=\mathcal{H}.
\end{eqnarray*}
Now we introduce the bilinear form given by
\begin{equation*}
	\mathfrak{a}\left(  (y, \alpha),(\varphi, \beta)\right):=\int_{0}^{L} y^{\prime}(x)\varphi^{\prime}(x) dx
\end{equation*}
on the domain $D(\mathfrak{a})=\mathcal{H}^{1}.$
\par 
In order to state the well-posedness of \eqref{s3}, we need some technical lemmas. Although based on elementary facts, the proofs of these lemmas are reported in the Appendix \ref{Appendix A}. The first of these lemmas is the following:
\begin{lemma} \label{Lemma 1}
	For all $t\in (0,T)$, the operator $\mathcal{B}(t)$ is uniformly bounded and its adjoint is given by 
	\begin{eqnarray*}
		(\mathcal{B}(t))^{*}(
		\varphi, \beta)=\left( -a(\cdot,t)\varphi
		, -\frac{\kappa}{\mu}\displaystyle\int_{0}^{L}b(x,t)\varphi(x)dx- c(t)\beta\right) \quad \forall (
		\varphi, \beta)\in\mathcal{H}.
	\end{eqnarray*}
\end{lemma}
The second auxiliary lemma is an important ingredient in the proof that $\mathcal{A}$ generates an analytic semigroup on $\mathcal{H}$.
\begin{lemma}\label{Lemma 2}
	The form $\mathfrak{a}$ is densely
	defined, closed symmetric and positive.
\end{lemma}
We are now in a position to prove that $\mathcal{A}$ is self-adjoint and generates an analytic semigroup on $\mathcal{H}$.
\begin{proposition}
	The operator $\mathcal{A}$ is densely defined, self-adjoint and generates an analytic $C_0$-semigroup $(e^{t\mathcal{A}})_{t\geq 0}$ on $\mathcal{H}$. Furthermore, we have the following interpolation result.  
	\begin{eqnarray*}
		\left[\mathcal{H},D(\mathcal{A})\right]_{1/2,2}=\mathcal{H}^{1}. \label{interpolation result}
	\end{eqnarray*}
\end{proposition}
\begin{proof}
	According to Lemma \ref{Lemma 2} and the second representation Theorem of \cite[Theorem 2.23]{kato2013perturbation}, the form $\mathfrak{a}$  induces a positive self-adjoint sectorial operator $\tilde{\mathcal{A}}$ on $\mathcal{H}$ which is given as follows. A couple $(y,\alpha)\in\mathcal{H}^{1}$ belongs to $D(\tilde{\mathcal{A}})$ if and
	only if there is $(f,c)\in \mathcal{H}$ such that 
	\begin{eqnarray}
		\mathfrak{a}\left( (y,\alpha),(\varphi,\beta)\right)=\left\langle (f,c),(\varphi,\beta)\right\rangle_{\mathcal{H}}, \qquad \forall (\varphi,\beta)\in \mathcal{H}^{1}. \label{R}
	\end{eqnarray}
	In this case $\tilde{\mathcal{A}}(y,\alpha)=(f,c)$. Furthermore, we have 
	\begin{eqnarray}
		D(\tilde{\mathcal{A}}^{1/2})=D(\mathfrak{a}). \label{Domain}
	\end{eqnarray}
	To conclude, we show that $\mathcal{A}=-\tilde{\mathcal{A}}$. Let $(y,\alpha)\in D(\tilde{\mathcal{A}})$ and choosing $\beta=0$ in \eqref{R}, we obtain 
	\begin{eqnarray*}
		\int_{0}^{L}y^{\prime}(x)\varphi^{\prime}(x)dx=\int_{0}^{L}f(x)\varphi(x) dx \quad \forall \varphi\in D(0,L).
	\end{eqnarray*}
	Consequently $y^{\prime\prime}\in L^{2}(0,L)$ and $-y^{\prime\prime}=f$ in $L^{2}(0,L)$. We return to \eqref{R}, by integration by parts and substituting the latter identity, we obtain
	\begin{eqnarray*}
		y^{\prime}(1)\varphi(1)-\mu y^{\prime}(0)\beta=\frac{\mu}{\kappa}c\beta\quad \forall  (\varphi,\beta)\in \mathcal{H}^{1}.
	\end{eqnarray*}
	This identity holds if and only if, $y^{\prime}(1)=0$ and $-\kappa y^{\prime}(0)=c$. In consequence, $D(\tilde{\mathcal{A}})\subset D(\mathcal{A})$ and $\mathcal{A}=-\tilde{\mathcal{A}}$. The other inclusion is by simple integration by parts. Hence, $-\mathcal{A}=\tilde{\mathcal{A}}$ and
	so $\mathcal{A}$ is self-adjoint and generates an analytic $C_0$-semigroup on $\mathcal{H}$. Finally, Theorem 4.36 of \cite{Lunardi} yields $D(\tilde{\mathcal{A}}^{1/2})=\left[\mathcal{H}, D(\tilde{\mathcal{A}})\right]_{1/2,2}$ which, combined with \eqref{Domain} and $\tilde{\mathcal{A}}=-\mathcal{A}$, shows the interpolation result in question. 
\end{proof}
\begin{remark} \label{Remark 1}
	Since $\mathcal{A}$ is self-adjoint, the homogeneous adjoint equation of \eqref{acp} is given by
	\begin{equation*}
		\left\{
		\begin{aligned}
			-\mathbf{U}^{\prime}(t)&=(\mathcal{A}+(\mathcal{B}(t))^{*}) \mathbf{U}(t), \quad 0<t<T,\\
			\mathbf{U}(T)&=\mathbf{U}_{T}.
		\end{aligned}
		\right.
	\end{equation*} 
	Consequently, the properties of $\mathcal{A}$ and $(\mathcal{B}(t))^{*}$, justify the form of the system \eqref{s2}.
\end{remark}
\par
We are now focusing on addressing the following solution categories for \eqref{s1}.
\begin{definition}[\textbf{Weak solution}]
	Let $(y_{0},z_{0})\in \mathcal{H}$ and $(f,g)\in L^{2}(0,T;\mathcal{H})$. A weak solution of \eqref{s1} is a couple of functions $(y,z)\in \mathcal{W}$ such that 
	\begin{eqnarray*}
		&&\int_{0}^{T}\langle y_{t},\varphi \rangle_{(H^{1}(0,L))^{\prime}, H^{1}(0,L)} \d t +\int_{Q}y_{x}\varphi_{x}\d x\d t+ \int_{Q}ay\varphi\d x\d t+ \int_{Q}bz\varphi\d x\d t \nonumber\\
		&&+ \mu\kappa^{-1}\int_{0}^{T}z^{\prime}\rho \d t +\mu\kappa^{-1}\int_{0}^{T}cz\rho \d t=\int_{Q}f\varphi \d x\d t +\mu\kappa^{-1}\int_{0}^{T}g\rho \d t\label{weak solution}
	\end{eqnarray*}
	for all $(\varphi,\rho)\in \mathcal{W}$ with $\varphi(\cdot,T)=\rho(T)=0$.
\end{definition}
We prove using \cite[Theorem 1.1, p. 37]{lions2013equations} the following existence results.
\begin{proposition}
	Let $(y_{0},z_{0})\in \mathcal{H}$ and $(f,g)\in L^{2}(0,T;\mathcal{H})$. The system \eqref{s1} has a unique weak solution $(y,z)\in \mathcal{W}$. Moreover, we have the estimate
	\begin{eqnarray*}
		&&\max_{0\leq t\leq T}\|(y(\cdot,t),z(t))\|^{2}_{\mathcal{H}}+ \|(y,z)\|^{2}_{L^{2}(0,T; \mathcal{H}^{1})}+ \|(y_{t},z^{\prime})\|^{2}_{L^{2}(0,T; (\mathcal{H}^{1})^{\prime})}\\
		&& \leq C\left( \|(y_{0},z_{0})\|^{2}_{\mathcal{H}}+ \|(f,g)\|^{2}_{L^{2}(0,T; \mathcal{H})}\right)
	\end{eqnarray*}
	for some positive constant $C$.
\end{proposition}
\begin{definition}[\textbf{Strong solution}]
	Let $(y_{0},z_{0})\in \mathcal{H}$ and $(f,g)\in L^{2}(0,T;\mathcal{H})$. 
	A strong solution of \eqref{s1} is a function $(y,z)\in \mathcal{E}$  fulfilling \eqref{s1} in $L^{2}(0,T;\mathcal{H})$.
\end{definition}

The following existence and uniqueness result is derived from \cite[Theorem 3.1]{pruss2001solvability}.
\begin{proposition}
	Let $(y_{0},z_{0})\in \mathcal{H}^{1}$ and $(f,g)\in L^{2}(0,T;\mathcal{H})$. The system \eqref{s1} has a unique strong solution $(y,z)\in \mathcal{E}$. Moreover, we have the estimate
	\begin{eqnarray*}
		&&\|(y,z)\|^{2}_{\mathcal{E}}\leq C\left( \|(y_{0},z_{0})\|^{2}_{\mathcal{H}}+ \|(f,g)\|^{2}_{L^{2}(0,T; \mathcal{H})}\right)
	\end{eqnarray*}
	for some positive constant $C$.
\end{proposition}
%\begin{proof}
%	It suffices to apply Theorem 3.1 in \cite{pruss2001solvability}.
%\end{proof}
\begin{remark} 
	It should be noted that we have not proved the existence and uniqueness of the system \eqref{s0}, nor the admissibility properties, as we have examined the intermediate system \eqref{s1} which allowed us to guarantee the existence and uniqueness of the system \eqref{s0}. Similar results can be established for the strong and weak solutions of the adjoint system \eqref{s2}.
	Furthermore, the duality relation between systems \eqref{s1} and \eqref{s2} is given by 
	\begin{eqnarray}
		\int_{\omega_{T}}v\varphi dxdt=\left\langle (
		y(\cdot,T), z(T))
		, (\varphi_{T},\rho_{T}) \right\rangle_{\mathcal{H}}-\left\langle (
		y_{0}, z_{0}), (\varphi(\cdot,0), \rho(0))\right\rangle_{\mathcal{H}} \label{dr}
	\end{eqnarray}
	for all weak solutions $(y,z)$ and $(\varphi,\rho)$ of \eqref{s1} and \eqref{s2} respectively. This relation will be shown in Appendix \ref{Appendix B}.
\end{remark}

\section{Carleman estimate and proof of main results} \label{sec3}
\subsection{Carleman estimate}
In order to study the observability of \eqref{s2}, we will establish a new Carleman estimate. First, let us recall the definitions of several classical weights, frequently used in this context. 
\par 
Let $\omega^{\prime}\subset\subset \omega$ be a nonempty open set and consider the following positive weight functions $\alpha, \hat{\alpha}, \xi$ and $\hat{\xi}$ which depend on $L$ and $\omega$
\begin{equation*}
	\begin{aligned}
		\alpha(x, t)&:=\frac{e^{2\lambda m}-e^{\lambda\left(m+ \eta(x)\right)}}{t(T-t)}\quad &&\forall (x,t)\in Q,\\
		\xi(x, t)&:=\frac{e^{\lambda\left(m+ \eta(x)\right)}}{t(T-t)}\quad\quad\quad &&\forall (x,t)\in Q,\\
		\hat{\alpha}(t)&:= \max_{x\in [0,L]}\alpha(x, t)=\alpha(0, t)=\alpha(L, t) \quad &&\forall t\in (0,T),\\
		\hat{\xi}(t)&:= \min_{x\in [0,L]}\xi(x, t)=\xi(0, t)=\xi(L, t) \quad &&\forall t\in (0,T).
	\end{aligned}
\end{equation*}
Here, $m>1$ and $\lambda>1$ is a sufficiently large positive constant (to be chosen later), and $\eta=\eta(x)$ is a function in $C^2([0,L])$ satisfying
\begin{equation}
	\label{eta}
	\eta>0 \text { in } (0,L), \quad \eta(0)=\eta(L)=0,\quad  \inf_{(0,L) \backslash \omega^{\prime}}\left|\eta^{\prime}\right| >0\quad\mbox{and}\quad \max_{[0,L]}\eta=1.
\end{equation}
The existence of such function $\eta$ satisfying \eqref{eta} is proved in \cite{fursikov1996controllability}. For the sake of brevity, the Lebesgue integration elements $\d x$ and $\d t$ will be omitted in Carleman's estimate, as will his proof. A global Carleman estimate holds for the solutions to a simplified version of \eqref{s2} with term sources: 
\begin{equation}\label{s4}
	\left\{
	\begin{aligned}
		&-\varphi_{t}-\varphi_{xx}=f
		& & \text {in}\;Q, \\
		& -\rho^{\prime}(t)-\kappa\varphi_{x}(0,t)
		=g & & \text {in}\;(0,T), \\
		& \varphi(0,t)=\mu \rho(t)  & & \text {in}\;(0,T), \\
		& \varphi_{x}(L,t) =0 & & \text {in}\;(0,T), \\
		& \varphi(\cdot,T)=\varphi_{T}  & & \text {in}\;(0,L), \\
		& \rho(T)=\rho_{T}   \\	
	\end{aligned}
	\right.
\end{equation}
can be stated as follows.
\begin{proposition} \label{P4}
	There are constants $C_{1}>0$ and $\lambda_1, s_{1}> 1$ depending only on $\omega$ and $L$ such that for any $s\geq s_{1}\left[\left(1+(\mu\kappa)^{-1}+(\mu\kappa)^{-1/2}\right)T+\left(1+\mu\kappa+(\mu\kappa)^{1/2}\right)T^{2}\right]$, any $\lambda\geq\lambda_1$ and any strong solution $(\varphi,\rho)$ of \eqref{s4}, we have the following estimate
	\begin{eqnarray}
		&&s^{-1}\int_{Q}e^{-2s\alpha}\xi^{-1}|\varphi_{t}|^{2} +s^{-1}\int_{Q}e^{-2s\alpha}\xi^{-1}|\varphi_{xx}|^{2}+ \mu\kappa^{-1}s^{-1}\int_{0}^{T}e^{-2s\hat{\alpha}}\hat{\xi}^{-1}|\rho^{\prime}|^{2} \nonumber\\
		&&+ s^{3} \lambda^{4} \int_{Q} e^{-2s\alpha}\xi^{3} |\varphi|^{2}  + s^{3}\lambda^{3}\int_{0}^{T}e^{-2s\hat{\alpha}}\hat{\xi}^{3}\left(|\varphi(0,t)|^{2}+ |\varphi(L,t)|^{2}\right) \nonumber\\
		&&+s \lambda^{2}\int_{Q}e^{-2s\alpha}\xi|\varphi_{x}|^{2}  + s\lambda\int_{0}^{T}e^{-s\hat{\alpha}}\hat{\xi}|\varphi_{x}(0,t)|^{2} \nonumber\\
		&&  \leq C\left( \int_{Q}e^{-2s\alpha}|f|^{2}+\mu\kappa^{-1}\int_{0}^{T}e^{-2s\hat{\alpha}}|g|^{2} +s^{3} \lambda^{4} \int_{\omega_T} e^{-s\alpha}\xi^{3} |\varphi|^{2}   \right). \label{Carleman estimate}
	\end{eqnarray}
\end{proposition}
This Carleman estimate is one of the main contributions of this paper and the key to proving its main results. The difficulty in this case is due to boundary and internal couplings that give unwanted terms; the choice of the inner product in \eqref{scalar product} will eliminate these terms. 
%The nonlocal term in space can be absorbed as a source term. 
%The proof of Proposition \ref{P4} is given in Appendix \ref{Appendix B}.
\begin{proof}
	We will denote by $C$ a generic positive constant that will be changed from one line to the next, $s_1>1$ and $\lambda_{1}>1$ are constants that will be increased from one passage to the next, the aim being to absorb terms. In the sequel $C, s_1$ and $\lambda_{1}$ will depend only on $\omega$ and $L$. For simplicity, the proof will be divided into several steps.\\
	\textbf{\underline{Step 1.} Change of unknowns.}
	Let $(\varphi,\rho)$ be a strong solution of \eqref{s4}. 
	Define 
	\begin{equation*}
		\begin{aligned}
			\psi&:= e^{-s \alpha}\varphi, \quad \tilde{f}:= -e^{-s \alpha}f\quad &&\forall (x,t)\in Q,\\
			\gamma&:= e^{-s \hat{\alpha}}\rho, \quad\; \tilde{g}:= -e^{-s \hat{\alpha}}g \quad &&\forall t\in (0,T).
		\end{aligned}
	\end{equation*}
	We have the following elementary identities:
	\begin{equation} \label{elementary identities}
		\begin{aligned}
			\alpha_{x} &=-\lambda\xi \eta^{\prime},\\
			\alpha_{xx} &= -\lambda^{2} \xi\left| \eta^{\prime}\right|^{2}-\lambda \xi \eta^{\prime\prime},\\
			\gamma^{\prime}(t) &= e^{-s\hat{\alpha}}\rho^{\prime}(t)-s\gamma(t)\hat{\alpha}_{t},\\
			\psi_{t} &= e^{-s\alpha}\varphi_{t}-s\psi\alpha_{t},\\
			\psi_{x}&= e^{-s \alpha}\varphi_{x}-s \psi \alpha_{x}=e^{-s \alpha}\varphi_{x}+s \lambda \xi\psi \eta^{\prime},\\
			\psi_{xx}&= e^{-s \alpha} \varphi_{xx}-2 s \psi_{x} \alpha_{x}-s^{2} \psi|\alpha_{x}|^{2}-s \psi \alpha_{xx}\\
			&= e^{-s \alpha} \varphi_{xx}+2s\lambda\xi\psi_{x}\eta^{\prime}-s^{2}\lambda^{2}\xi^{2}\psi|\eta^{\prime}|^{2}+s\lambda^{2}\xi\psi|\eta^{\prime}|^{2}+s\lambda\xi\psi\eta^{\prime\prime}.
		\end{aligned}
	\end{equation}
	The heat equation is transformed on $Q$ as follows 
	\begin{eqnarray*}
		\psi_{t}+\psi_{xx}=\tilde{f}-s\psi\alpha_{t} +2s\lambda\xi\psi_{x}\eta^{\prime}-s^{2}\lambda^{2}\xi^{2}\psi|\eta^{\prime}|^{2}+s\lambda^{2}\xi\psi|\eta^{\prime}|^{2}+s\lambda\xi\psi\eta^{\prime\prime},
	\end{eqnarray*}
	while the ODE on $(0,T)$ yields, 
	\begin{eqnarray*}
		\gamma^{\prime}(t)+\kappa\psi_{x}(0,t) =\tilde{g}-s\gamma(t)\hat{\alpha}_{t}+\kappa s\lambda\eta^{\prime}(0)\hat{\xi}\psi(0,t).
	\end{eqnarray*}
	As usual, we take the following decomposition:
	\begin{equation*}
		\begin{aligned}
			&M_{1}\psi+ M_{2}\psi=\tilde{f}+s\lambda\psi\xi\eta^{\prime\prime} -s\lambda^{2}\psi\xi|\eta^{\prime}|^{2}:=F\;&&\text{in}\;Q,\\
			&N_{1}(\psi,\gamma)+ N_{2}(\psi,\gamma)=\tilde{g} \;&&\text{in}\;(0,T),
		\end{aligned}
	\end{equation*}
	where
	\begin{equation} \label{m1m2n1}
		\begin{aligned}
			M_{1}\psi&:=-2s\lambda^{2}\psi\xi|\eta^{\prime}|^{2}-2s\lambda\xi\psi_{x}\eta^{\prime}+\psi_{t}:=\sum_{j=1}^{3}(M_{1}\psi)_{j},\\
			M_{2}\psi &:= s^{2}\lambda^{2}\psi\xi^{2}|\eta^{\prime}|^{2}+\psi_{xx} +s\psi\alpha_{t}:=\sum_{j=1}^{3}(M_{2}\psi)_{j},\\
			%						F&:=\tilde{f}+s\lambda\psi\xi\eta^{\prime\prime} -s\lambda^{2}\psi\xi|\eta^{\prime}|^{2},\\
			N_{1}(\psi,\gamma) &:=\gamma^{\prime}(t)-\kappa s\lambda\eta^{\prime}(0)\hat{\xi}\psi(0,t):=\sum_{j=1}^{2}(N_{1}(\psi,\gamma))_{j},	\\
			N_{2}(\psi,\gamma)&:= s\gamma(t)\hat{\alpha}_{t}+\kappa\psi_{x}(0,t):=\sum_{j=1}^{2}(N_{2}(\psi,\gamma))_{j}.
			%						G&:=\tilde{g}+b(t)\gamma(t)+\textcolor{red}{\frac{\kappa}{\mu}}\int_{0}^{L}a(x,t)\psi(x,t)dx.
		\end{aligned}
	\end{equation}
	Let us apply the parallelogram identity in $L^{2}(Q)$ and $L^{2}(0,T)$, we obtain
	\begin{eqnarray}
		\|F\|^{2}_{L^{2}(Q)}= \|M_{1}\psi\|^{2}_{L^{2}(Q)}+\|M_{2}\psi\|^{2}_{L^{2}(Q)}+ 2\sum_{1\leq j, k\leq 3 }\langle (M_{1}\psi)_{j}, (M_{2}\psi)_{k}\rangle_{L^{2}(Q)} \label{c1}
	\end{eqnarray}
	and
	\begin{eqnarray}
		\|\tilde{g}\|^{2}_{L^{2}(0,T)}&=& \|N_{1}(\psi,\gamma)\|^{2}_{L^{2}(0,T)}+ \|N_{2}(\psi,\gamma)\|^{2}_{L^{2}(0,T)}\nonumber\\
		&& + 2\sum_{1\leq j, k\leq 2}\langle (N_{1}(\psi,\gamma))_{j}, (N_{2}(\psi,\gamma))_{k}\rangle_{L^{2}(0,T)}. \label{c2}
	\end{eqnarray}
	Furthermore, the third equation of \eqref{s0} yields
	\begin{eqnarray}
		\psi(0,t)=\mu\gamma(t)\quad\mbox{and}\quad
		\psi_{x}(L,t)=s \lambda \eta^{\prime}(L)\hat{\xi}(t)\psi(L,t) . \label{FC}
	\end{eqnarray}
	\textbf{\underline{Step 2.} Estimating the mixed terms in \eqref{c1} and \eqref{c2} from below.} We often use the following basic pointwise estimates on $\overline{Q}$,
	\begin{equation*}
		|\alpha_{x}| \leq C \lambda \xi, \quad\left| \alpha_{t}\right| \leq C T\xi^{2}, \quad \left| \alpha_{tt}\right| \leq C T^{2}\xi^{3}, \quad \left| \xi_{t}\right| \leq C T\xi^{2}. \label{basic estimates}
	\end{equation*}
	\textbf{\underline{Step 2a.} Estimate from below of $\langle M_{1}\psi, (M_{2}\psi)_{1}\rangle_{L^{2}(Q)}$.}
	It is obvious that 
	\begin{eqnarray*}
		\langle M_{1}\psi, (M_{2}\psi)_{1}\rangle_{L^{2}(Q)}=\sum_{1\leq j \leq 3} \langle (M_{1}\psi)_{j}, (M_{2}\psi)_{1}\rangle_{L^{2}(Q)}.
	\end{eqnarray*}
	Firstly, we have
	\begin{eqnarray*}
		\langle (M_{1}\psi)_{1}, (M_{2}\psi)_{1}\rangle_{L^{2}(Q)}=-2s^{3}\lambda^{4}\int_{Q}|\eta^{\prime}|^{4}\xi^{3}\psi^{2}.
	\end{eqnarray*}
	Using integration by parts and $\xi_{x}=\lambda\xi\eta^{\prime}$, we further derive
	\begin{eqnarray*}
		\left\langle\left(M_{1} \psi\right)_{2},\left(M_{2} \psi\right)_{1}\right\rangle_{L^{2}\left(Q\right)}
		%					&=&-2s^{3}\lambda^{3}\int_{Q}(\eta^{\prime})^{3}\xi^{3}\psi \psi_{x}\\
		%					&=&-s^{3} \lambda^{3} \int_{Q}\left(\eta^{\prime}\right)^{3} \xi^{3}   \left(\psi^{2}\right)_{x}  \\
		%					&=&s^{3} \lambda^{3} \int_{Q} \left(\left(\eta^{\prime}\right)^{3} \xi^{3} \right)_{x} \psi^{2} -s^{3} \lambda^{3} \int_{0}^{T}\left[\left(\eta^{\prime}\right)^{3} \xi^{3} \psi^{2} \right]_{x=0}^{x=L} \\
		&=& 3s^{3} \lambda^{4} \int_{Q}\left|\eta^{\prime}\right|^{4}\xi^{3} \psi^{2} + 3s^{3} \lambda^{3} \int_{Q} \eta^{\prime\prime}\left|\eta^{\prime}\right|^{2}\xi^{3} \psi^{2}\\
		&&-s^{3} \lambda^{3} \int_{0}^{T}\left[\left(\eta^{\prime}\right)^{3} \xi^{3} \psi^{2} \right]_{x=0}^{x=L}.
	\end{eqnarray*}
	From the properties of $\eta$, it follows that 
	$$\int_{Q}\left|\eta^{\prime}\right|^{4} \xi^{3} \psi^{2}\geqslant \int_{Q\backslash \omega^{\prime}_T}\left|\eta^{\prime}\right|^{4} \xi^{3} \psi^{2}\geqslant C\int_{Q}\xi^{3} \psi^{2}-C\int_{ \omega^{\prime}_T}\xi^{3} \psi^{2}.$$
	For sufficiently large $\lambda_{1}$, we obtain
	\begin{eqnarray}
		&&\left\langle\left(M_{1} \psi\right)_{1},\right. \left.\left(M_{2} \psi\right)_{1}\right\rangle_{L^{2}\left(Q\right)}+\left\langle\left(M_{1} \psi\right)_{2},\left(M_{2} \psi\right)_{1}\right\rangle_{L^{2}\left(Q\right)} \geq  Cs^{3} \lambda^{4} \int_{Q} \xi^{3} \psi^{2} \nonumber\\
		&& \quad\quad -Cs^{3} \lambda^{4} \int_{\omega^{\prime}_T} \xi^{3} \psi^{2} -s^{3} \lambda^{3} \int_{0}^{T}\left[\left(\eta^{\prime}\right)^{3} \xi^{3} \psi^{2} \right]_{x=0}^{x=L}.
		\label{c3}
	\end{eqnarray}
	Integrating by parts in time, with \eqref{c2}, $\psi(\cdot,0)=\psi(\cdot,T)=0$ and $\left| \xi_{t}\right| \leq C T\xi^{2}$, yields
	\begin{eqnarray*}
		\left\langle\left(M_{1} \psi\right)_{3},\left(M_{2} \psi\right)_{1}\right\rangle_{L^{2}\left(Q\right)}
		%					&=&\frac{1}{2}s^{2} \lambda^{2} \int_{Q}\left|\eta^{\prime}\right|^{2} \xi^{2} \left(\psi^{2}\right)_{t} \\
		=-s^{2} \lambda^{2} \int_{Q}\left|\eta^{\prime}\right|^{2}  \xi_{t} \xi \psi^{2} \geq -C Ts^{2} \lambda^{2} \int_{Q} \xi^{3} \psi^{2}.
	\end{eqnarray*}
	This term is absorbed by the first term in the right-hand side of \eqref{c3} if we take $s\geq s_1 T$ and $\lambda\geq\lambda_{1}>1$. Altogether, we have shown
	\begin{eqnarray}
		\sum_{1\leqslant j\leqslant3}\left\langle\left(M_{1} \psi\right)_{j},\left(M_{2} \psi\right)_{1}\right\rangle_{L^{2}\left(Q\right)}
		& \geq& C s^{3} \lambda^{4} \int_{Q} \xi^{3} \psi^{2}-C s^{3} \lambda^{4} \int_{\omega^{\prime}_T} \xi^{3} \psi^{2} \nonumber\\
		&&-s^{3} \lambda^{3} \int_{0}^{T}\left[\left(\eta^{\prime}\right)^{3} \xi^{3} \psi^{2} \right]_{x=0}^{x=L}, \label{step2a}
	\end{eqnarray}
	for any $s\geq s_{1}T$ and any $\lambda\geq\lambda_{1}$.\\
	\textbf{\underline{Step 2b.} Estimate from below of $\langle M_{1}\psi, (M_{2}\psi)_{2}\rangle_{L^{2}(Q)}$.}
	It is obvious that 
	\begin{eqnarray*}
		\langle M_{1}\psi, (M_{2}\psi)_{2}\rangle_{L^{2}(Q)}=\sum_{1\leq j \leq 3} \langle (M_{1}\psi)_{j}, (M_{2}\psi)_{2}\rangle_{L^{2}(Q)}.
	\end{eqnarray*}
	By integration by parts and $\xi_{x}=\lambda\xi\eta^{\prime}$, we obtain
	\begin{eqnarray*}
		\left\langle\left(M_{1} \psi\right)_{1},\right.  \left.\left(M_{2} \psi\right)_{2}\right\rangle_{L^{2}\left(Q\right)}
		%					&=&-2 s \lambda^{2} \int_{Q}\left|\eta^{\prime}\right|^{2} \xi \psi \psi_{xx} \\
		%					&=& 2s \lambda^{2} \int_{Q} \left(\left| \eta^{\prime}\right|^{2} \xi \psi\right)_{x}\psi_{x}-2 s \lambda^{2} \int_{0}^{T}\left[\left|\eta^{\prime}\right|^{2} \xi \psi \psi_{x}\right]_{x=0}^{x=L} \\
		&=&2s \lambda^{2} \int_{Q}  \left(\left| \eta^{\prime}\right|^{2}\right)_{x}\xi \psi\psi_{x} +2 s \lambda^{3} \int_{Q}\left( \eta^{\prime}\right)^{3}\xi \psi  \psi_{x} \\
		&& +2s \lambda^{2} \int_{Q}\left| \eta^{\prime}\right|^{2} \xi|\psi_{x}|^{2} -2 s \lambda^{2} \int_{0}^{T}\left[\left|\eta^{\prime}\right|^{2} \xi \psi \psi_{x}\right]_{x=0}^{x=L}.
	\end{eqnarray*}
	As above, the third summand will lead to a term controlling $| \psi_{x}|^{2}$. We now apply Young's inequality to the first and second terms, $\xi\geq \frac{C}{T^{2}}$, we find 
	\begin{eqnarray*}
		\left|2 s \lambda^{2} \int_{Q}  \left(\left| \eta^{\prime}\right|^{2}\right)_{x}\xi \psi\psi_{x}\right| &\leq& C s \lambda^{4} \int_{Q} \xi |\psi|^{2}+C s \int_{Q} \xi|\psi_{x}|^{2}\\
		&\leq& C s^{2} \lambda^{4} \int_{Q} \xi^{2} |\psi|^{2}+C s \int_{Q} \xi|\psi_{x}|^{2},
	\end{eqnarray*}
	for any $s\geq s_1 T^{2}$, and
	\begin{eqnarray*}
		\left|2 s \lambda^{3} \int_{Q}\left( \eta^{\prime}\right)^{3}\xi \psi  \psi_{x}\right|&\leq& C  s^{2} \lambda^{4} \int_{Q} \xi^{2} |\psi|^{2}+C  \lambda^{2} \int_{Q}\left|\psi_{x}\right|^{2}.
	\end{eqnarray*}
	It follows
	\begin{eqnarray*}
		\left\langle\left(M_{1} \psi\right)_{1},\left(M_{2} \psi\right)_{2}\right\rangle_{L^{2}\left(Q\right)} &\geq& -Cs^{2} \lambda^{4}\int_{Q}\xi^{3} |\psi|^{2}-C\int_{Q}\left(s \xi+\lambda^{2}\right)|\psi_{x}|^{2},\\
		&& +2s \lambda^{2} \int_{Q}\left| \eta^{\prime}\right|^{2} \xi|\psi_{x}|^{2} -2 s \lambda^{2} \int_{0}^{T}\left[\left|\eta^{\prime}\right|^{2} \xi \psi \psi_{x}\right]_{x=0}^{x=L}.
	\end{eqnarray*}
	Using $\xi_{x}=\lambda\xi\eta^{\prime}$, the next summand is given by
	\begin{eqnarray*}
		\left\langle\left(M_{1} \psi\right)_{2},\left(M_{2} \psi\right)_{2}\right\rangle_{L^{2}\left(Q\right)}
		%					&=&-2s \lambda \int_{Q} \eta^{\prime}\xi\psi_{x} \psi_{xx}=-s \lambda \int_{Q} \eta^{\prime}\xi(|\psi_{x}|^{2})_{x} \\
		%					&=&- s \lambda\int_{0}^{T} \left[\eta^{\prime}\xi |\psi_{x}|^{2}\right]_{x=0}^{x=L} +  s \lambda\int_{Q}(\eta^{\prime}\xi)_{x}|\psi_{x}|^{2} \\
		&=&- s \lambda\int_{0}^{T} \left[\eta^{\prime}\xi |\psi_{x}|^{2}\right]_{x=0}^{x=L}+  s \lambda\int_{Q}\eta^{\prime\prime}\xi|\psi_{x}|^{2} \\
		&&+s \lambda^{2}\int_{Q}|\eta^{\prime}|^{2}\xi|\psi_{x}|^{2}.
	\end{eqnarray*}
	Integrating by parts in space and $\psi_{x}(\cdot,0)=\psi_{x}(\cdot,T)=0$, we obtain
	\begin{eqnarray*}
		\left\langle\left(M_{1} \psi\right)_{3},\left(M_{2} \psi\right)_{2}\right\rangle_{L^{2}\left(Q\right)}
		%					&=&\int_{Q}  \psi_{t}\psi_{xx}=\int_{0}^{T}  \left[\psi_{t}\psi_{x}\right]_{x=0}^{x=L}-\int_{Q}  \psi_{xt}\psi_{x}\\
		&=& \int_{0}^{T}  \left[\psi_{t}\psi_{x}\right]_{x=0}^{x=L}-\frac{1}{2}\int_{Q}  (|\psi_{x}|^{2})_{t}\\
		&=& \int_{0}^{T}  \left[\psi_{t}\psi_{x}\right]_{x=0}^{x=L}.
	\end{eqnarray*}
	Using the fact that $|\eta^{\prime}|\geq C$ in $\overline{(0,L)\backslash\omega^{\prime}}$,  %$|\nabla\eta^0|=|\partial_{\nu}\eta^0|$ and $|\nabla\psi|^{2}=|\nabla_{\Gamma}\psi|^{2}+|\partial_{\nu}\psi|^{2}$ on $\Gamma$,
	we arrive at:
	\begin{eqnarray}
		&&\sum_{1\leqslant j\leqslant3}\left\langle\left(M_{1} \psi\right)_{j},\left(M_{2} \psi\right)_{2}\right\rangle_{L^{2}\left(Q\right)}
		\geq Cs \lambda^{2}\int_{Q}\xi|\psi_{x}|^{2}-C s^{2} \lambda^{4}\int_{Q}\xi^{3} \psi^{2} \nonumber\\
		&&-C \int_{Q}\left(s\xi+\lambda^{2}\right)|\psi_{x}|^{2}-Cs \lambda^{2}\int_{\omega^{\prime}_{T}}\xi|\psi_{x}|^{2} -2 s \lambda^{2} \int_{0}^{T}\left[\left|\eta^{\prime}\right|^{2} \xi \psi \psi_{x}\right]_{x=0}^{x=L} \nonumber\\
		&& - s \lambda\int_{0}^{T} \left[\eta^{\prime}\xi |\psi_{x}|^{2}\right]_{x=0}^{x=L} +\int_{0}^{T}  \left[\psi_{t}\psi_{x}\right]_{x=0}^{x=L}, \label{step2b}
	\end{eqnarray}
	for any $s\geq s_1 T^{2}$ and any $\lambda\geq \lambda_1$.\\
	\textbf{\underline{Step 2c.} Estimate from below of $\langle M_{1}\psi, (M_{2}\psi)_{3}\rangle_{L^{2}(Q)}$.}
	From $|\alpha_{t}|\leqslant CT\xi^2$, we get
	\begin{eqnarray*}
		\left\langle\left(M_{1} \psi\right)_{1},\left(M_{2} \psi\right)_{3}\right\rangle_{L^{2}\left(Q\right)} =-2 s^{2} \lambda^{2} \int_{Q}\left|\eta^{\prime}\right|^{2}  \alpha_{t} \xi \psi^{2} 
		\geq-Cs^{3} \lambda^{2} \int_{Q} \xi^{3} \psi^{2},
	\end{eqnarray*}
	for any $s\geq s_{1}T$. Using also $| \alpha_{xt}|\leqslant CT\lambda\xi^2$, $\xi_{x}=\lambda\xi\eta^{\prime}$ and integration by parts, we obtain
	\begin{eqnarray*}
		\left\langle\left(M_{1} \psi\right)_{2},\left(M_{2} \psi\right)_{3}\right\rangle_{L^{2}\left(Q\right)}
		%					&=&-2s^2 \lambda \int_{Q}  \alpha_{t}\eta^{\prime} \xi \psi\psi_{x}  =-s^{2} \lambda \int_{Q}  \alpha_{t}\eta^{\prime} \xi\left(\psi^{2}\right)_{x}  \\
		%					&=&-s^{2} \lambda \int_{0}^{T}  \left[\alpha_{t}\eta^{\prime} \xi\psi^{2}\right]_{x=0}^{x=L} +s^{2} \lambda \int_{Q}\left( \alpha_{t}\eta^{\prime} \xi\right)_{x} \psi^{2} \\
		&= & -s^{2} \lambda \int_{0}^{T}  \left[\alpha_{t}\eta^{\prime} \xi\psi^{2}\right]_{x=0}^{x=L}+ s^{2} \lambda \int_{Q}\alpha_{xt}\eta^{\prime} \xi \psi^{2}  \\
		&& +s^{2} \lambda \int_{Q}\alpha_{t}\eta^{\prime\prime} \xi \psi^{2} +s^{2} \lambda^{2} \int_{Q}\alpha_{t}|\eta^{\prime}|^{2} \xi \psi^{2}  \\
		&\geq& -s^{2} \lambda \int_{0}^{T}  \left[\alpha_{t}\eta^{\prime} \xi\psi^{2}\right]_{x=0}^{x=L}-Cs^{3} \lambda^{2} \int_{Q} \xi^{3} \psi^{2},
	\end{eqnarray*}
	for any $s\geq s_{1}T$ and any $\lambda\geq\lambda_{1}$. Integrating by parts with respect to time, $\psi(\cdot,0)=\psi(\cdot,T)=0$ and $\left|\alpha_{tt}\right| \leq C T^{2}\xi^{3}$, we can derive
	\begin{eqnarray*}
		\left\langle\left(M_{1} \psi\right)_{3},\left(M_{2} \psi\right)_{3}\right\rangle_{L^{2}\left(Q\right)}
		%					 &=&s\int_{Q} \alpha_{t}\psi_{t}\psi dxdt=\frac{s}{2} \int_{Q} \alpha_{t} \left(\psi^{2}\right)_{t} 
		=-\frac{s}{2} \int_{Q} \alpha_{tt} \psi^{2}  
		\geq-C s^{2} \int_{Q} \xi^{3} \psi^{2},
	\end{eqnarray*}
	for any $s\geq s_1 T^{2}$.
	We conclude from the above inequalities that
	\begin{eqnarray}
		\sum_{1\leqslant j\leqslant3}\left\langle\left(M_{1} \psi\right)_{j},\left(M_{2} \psi\right)_{3}\right\rangle_{L^{2}\left(Q\right)} &\geq& -Cs^{3} \lambda^{2} \int_{Q} \xi^{3} \psi^{2}-s^{2} \lambda \int_{0}^{T}  \left[\alpha_{t}\eta^{\prime} \xi\psi^{2}\right]_{x=0}^{x=L} \nonumber\\
		&&-Cs^{2} \lambda^{2} \int_{Q} \xi^{3} \psi^{2}, \label{step2c}
	\end{eqnarray}
	for any $s\geq s_{1}(T+T^{2})$ and $\lambda\geq\lambda_{1}$.\\
	\textbf{\underline{Step 2d.} Estimate from below of $\langle N_{1}(\psi,\gamma), N_{2}(\psi,\gamma)\rangle_{L^{2}(0,T)}$.} We now consider the boundary terms $N_{1}$ and $N_{2}$, using an integration by parts in time, $\gamma(0)=\gamma(T)=0$ and $\psi(0,t)=\mu\gamma$. Using $\left|\hat{\alpha}_{tt}\right| \leq C T^{2}\hat{\xi}^{3}$, we obtain
	\begin{eqnarray*}	 
		\left\langle\left(N_{1}(\psi, \gamma)\right)_{1},\left(N_{2}\left(\psi, \gamma\right)\right)_{1}\right\rangle_{L^{2}\left(0,T\right)}&=&s\int_{0}^{T}\hat{\alpha}_{t}\gamma\gamma^{\prime}=-\frac{s}{2}\int_{0}^{T}\hat{\alpha}_{tt}|\gamma|^{2}\\
		&\geq& -C\mu^{-2}sT^{2}\int_{0}^{T}\hat{\xi}^{3}|\psi(0,t)|^{2}.
	\end{eqnarray*}
	Next, 
	\begin{eqnarray*}	 
		\left\langle\left(N_{1}(\psi, \gamma)\right)_{1},\left(N_{2}\left(\psi, \gamma\right)\right)_{2}\right\rangle_{L^{2}\left(0,T\right)}&=&\kappa\int_{0}^{T}\gamma^{\prime}\psi_{x}(0,t)=\mu^{-1}\kappa\int_{0}^{T}\psi_{t}(0,t)\psi_{x}(0,t).
	\end{eqnarray*}		
	Using $\left| \hat{\alpha}_{t}\right| \leq C T\hat{\xi}^{2}$, we obtain 
	\begin{eqnarray*}
		\left\langle\left(N_{1}(\psi, \gamma)\right)_{2},\left(N_{2}\left(\psi, \gamma\right)\right)_{1}\right\rangle_{L^{2}\left(0,T\right)}&=&-\kappa s^{2}\lambda\eta^{\prime}(0)\int_{0}^{T}\hat{\alpha}_{t}\hat{\xi}\psi(0,t)\gamma \\
		%					&=& -\textcolor{red}{\mu^{-1}\kappa}s^{2}\lambda\eta^{\prime}(0)\int_{0}^{T}\hat{\alpha}_{t}\hat{\xi}|\psi(0,t)|^{2}\\
		&\geq & -C \mu^{-1}\kappa s^{3}\lambda\int_{0}^{T}\hat{\xi}^{3}|\psi(0,t)|^{2},
	\end{eqnarray*}
	for any $s\geq s_{1}T$.
	Finally,
	\begin{eqnarray*}
		\left\langle\left(N_{1}(\psi, \gamma)\right)_{2},\left(N_{2}\left(\psi, \gamma\right)\right)_{2}\right\rangle_{L^{2}\left(0,T\right)}=-\kappa^{2} s \lambda\eta^{\prime}(0)\int_{0}^{T} \hat{\xi}\psi(0,t)\psi_{x}(0,t).
	\end{eqnarray*}
	As a result 
	\begin{eqnarray}
		&&\left\langle N_{1}(\psi, \gamma),N_{2}\left(\psi, \gamma\right)\right\rangle_{L^{2}\left(0,T\right)}\geq -C\mu^{-2}sT^{2}\int_{0}^{T}\hat{\xi}^{3}|\psi(0,t)|^{2} \nonumber\\
		&&+\mu^{-1}\kappa\int_{0}^{T}\psi_{t}(0,t)\psi_{x}(0,t)-C \mu^{-1}\kappa s^{3}\lambda\int_{0}^{T}\hat{\xi}^{3}|\psi(0,t)|^{2}\nonumber\\
		&&-\kappa^{2} s \lambda\eta^{\prime}(0)\int_{0}^{T} \hat{\xi}\psi(0,t)\psi_{x}(0,t),  \label{step2d}
	\end{eqnarray}
	for any $s\geq s_{1}(T+T^{2})$.\\
	\textbf{\underline{Step 3.} The transformed estimate.} \\
	Taking in account \eqref{step2a}, \eqref{step2b}, \eqref{step2c} and \eqref{step2d}, we obtain
	\begin{eqnarray}
		&&\langle M_{1}\psi, M_{2}\psi\rangle_{L^{2}(Q)}+\frac{\mu}{\kappa}\left\langle N_{1}(\psi, \gamma),N_{2}\left(\psi, \gamma\right)\right\rangle_{L^{2}\left(0,T\right)} \geq C s^{3} \lambda^{4} \int_{Q} \xi^{3} \psi^{2} \nonumber\\
		&& -C s^{3} \lambda^{4} \int_{\omega^{\prime}_T} \xi^{3} \psi^{2}-s^{3} \lambda^{3} \int_{0}^{T}\left[\left(\eta^{\prime}\right)^{3} \xi^{3} \psi^{2} \right]_{x=0}^{x=L} +Cs \lambda^{2}\int_{Q}\xi|\psi_{x}|^{2} \nonumber\\
		&&-C s^{2} \lambda^{4}\int_{Q}\xi^{3} \psi^{2}-C \int_{Q}\left(s\xi+\lambda^{2}\right)|\psi_{x}|^{2} -Cs \lambda^{2}\int_{\omega^{\prime}_{T}}\xi|\psi_{x}|^{2} \nonumber\\
		&&-2 s \lambda^{2} \int_{0}^{T}\left[\left|\eta^{\prime}\right|^{2} \xi \psi \psi_{x}\right]_{x=0}^{x=L} - s \lambda\int_{0}^{T} \left[\eta^{\prime}\xi |\psi_{x}|^{2}\right]_{x=0}^{x=L} +\int_{0}^{T}  \left[\psi_{t}\psi_{x}\right]_{x=0}^{x=L} \nonumber\\
		&& -Cs^{3} \lambda^{2} \int_{Q} \xi^{3} \psi^{2}-s^{2} \lambda \int_{0}^{T}  \left[\alpha_{t}\eta^{\prime} \xi\psi^{2}\right]_{x=0}^{x=L} -Cs^{2} \lambda^{2} \int_{Q} \xi^{3} \psi^{2} \nonumber\\
		&&-C\mu^{-1}\kappa^{-1}sT^{2}\int_{0}^{T}\hat{\xi}^{3}|\psi(0,t)|^{2}+\int_{0}^{T}\psi_{t}(0,t)\psi_{x}(0,t) \nonumber\\
		&&-C s^{3}\lambda\int_{0}^{T}\hat{\xi}^{3}|\psi(0,t)|^{2}-\mu\kappa s \lambda\eta^{\prime}(0)\int_{0}^{T} \hat{\xi}\psi(0,t)\psi_{x}(0,t). \label{trans estimate}
	\end{eqnarray}
	Using \eqref{FC} and $|\hat{\xi}|\leq CT\hat{\xi}^{2}$, we obtain
	\begin{eqnarray}
		\int_{0}^{T}  \left[\psi_{t}\psi_{x}\right]_{x=0}^{x=L}+\int_{0}^{T}\psi_{t}(0,t)\psi_{x}(0,t)&=&\int_{0}^{T}\psi_{t}(L,t)\psi_{x}(L,t) \nonumber\\
		&=& s \lambda \eta^{\prime}(L)\int_{0}^{T}\hat{\xi}(t)\psi_{t}(L,t)\psi(L,t) \nonumber\\
		&=& -\frac{s \lambda \eta^{\prime}(L)}{2}\int_{0}^{T}\hat{\xi}_{t}(t)\psi^{2}(L,t)  \nonumber\\
		&\geq& -Cs^{2}\lambda\int_{0}^{T}\hat{\xi}^{2}(t)\psi^{2}(L,t), \label{trans1}
	\end{eqnarray}
	for any $s\geq s_{1}T$. Using the estimates \eqref{trans estimate}, \eqref{trans1} and the fact that $\eta^{\prime}(0)<0,\; \eta^{\prime}(1)>0$, $|\eta^{\prime}|\leq C$, $|\hat{\alpha}_{t}|\leq CT\hat{\xi}^{2}$,  we deduce
	\begin{eqnarray}
		&&\|M_{1}\psi\|^{2}_{L^{2}(Q)}+\|M_{2}\psi\|^{2}_{L^{2}(Q)}+ \mu\kappa^{-1}\|N_{1}(\psi,\gamma)\|^{2}_{L^{2}(0,T)}+ \mu\kappa^{-1}\|N_{2}(\psi,\gamma)\|^{2}_{L^{2}(0,T)} \nonumber\\
		&&+ s^{3} \lambda^{4} \int_{Q} \xi^{3} |\psi|^{2} + s^{3}\lambda^{3}\int_{0}^{T}\hat{\xi}^{3}(t)\left(|\psi(0,t)|^{2}+ |\psi(L,t)|^{2}\right) \nonumber\\
		&&+s \lambda^{2}\int_{Q}\xi|\psi_{x}|^{2} + s\lambda\int_{0}^{T}\hat{\xi}(t)\left(|\psi_{x}(0,t)|^{2}+ |\psi_{x}(L,t)|^{2}\right) \nonumber\\
		&&  \leq C\left( \|F\|^{2}_{L^{2}(Q)}+ \mu\kappa^{-1}\|\tilde{g}\|^{2}_{L^{2}(0,T)} +s^{3} \lambda^{4} \int_{\omega^{\prime}_T} \xi^{3} |\psi|^{2} + s^{2} \lambda^{4}\int_{Q}\xi^{2} |\psi|^{2} \right. \nonumber\\
		&& \left. + s \lambda^{2}\int_{\omega^{\prime}_{T}}\xi|\psi_{x}|^{2} + s \lambda^{2} \int_{0}^{T}\hat{\xi}(t)|\psi(L,t)||\psi_{x}(L,t)|  + s^{2}\lambda\int_{0}^{T}\hat{\xi}^{2}(t)|\psi(L,t)|^{2} \right. \nonumber\\
		&& \left.  + s^{3}\lambda^{2}\int_{Q}\xi^{3}|\psi|^{2} + s^{2} \lambda T\int_{0}^{T}  \hat{\xi}^{3}(t)\left( |\psi(L,t)|^{2}+|\psi(0,t)|^{2}\right) +s^{2} \lambda^{2} \int_{Q} \xi^{3} |\psi|^{2}\right. \nonumber\\
		&& \left.   + \mu^{-1}\kappa^{-1}sT^{2}\int_{0}^{T}\hat{\xi}^{3}(t)|\psi(0,t)|^{2} +s^{3}\lambda \int_{0}^{T} \hat{\xi}^{3}(t)|\psi(0,t)|^{2} \right. \nonumber\\
		&& \left. +(\mu\kappa s\lambda+s\lambda^{2}) \int_{0}^{T} \hat{\xi}(t)|\psi(0,t)| |\psi_{x}(0,t)| + \int_{Q}\left(s\xi+\lambda^{2}\right)|\psi_{x}|^{2} \right). \label{trans2}
	\end{eqnarray}
	On the other hand, for the source term $F$, we have 
	\begin{eqnarray*}
		\|F\|^{2}_{L^{2}(Q)}&\leq & C\left(\int_{Q}e^{-2s\alpha}|f|^{2}+s^{2}\lambda^{4}\int_{Q}\xi^{2}|\psi|^{2}\right). \label{Source F}
	\end{eqnarray*}
	%				For the source term $G$, using the Cauchy Schwarz inequality in space for the nonlocal term and $\xi\geq \frac{C}{T^{2}}$, we obtain
	%				\begin{eqnarray}
		%					\|G\|^{2}_{L^{2}(0,T)}&\leq & C\left(\int_{0}^{T}e^{-2s\hat{\alpha}}|g|^{2}+\kappa^{-2}\|b\|^{2}_{L^{\infty}(0,T)}\int_{0}^{T}|\psi(0,t)|^{2}+\textcolor{red}{\mu^{-2}\kappa^{2}} \|a\|^{2}_{L^{\infty}(Q)}\int_{Q}\psi^{2}\right) \nonumber\\
		%					&\leq & C\left(\int_{0}^{T}e^{-2s\hat{\alpha}}|g|^{2}dt+\kappa^{-2}\|b\|^{2}_{L^{\infty}(0,T)}T^{6}\int_{0}^{T}\hat{\xi}^{3}(t)|\psi(0,t)|^{2}dt \right. \nonumber\\
		%					&& \left. + \mu^{-2}\kappa^{2}\|a\|^{2}_{L^{\infty}(Q)}T^{6}\int_{Q}\xi^{3}\psi^{2}dxdt\right). \label{Source G}
		%				\end{eqnarray}
	%				Multiplying \eqref{Source G} by $\textcolor{red}{\mu\kappa^{-1}}$, then taking $s\geq s_{1}\left(\mu^{-1/3}\kappa^{1/3}\|a\|^{2/3}_{L^{\infty}(Q)}+\mu^{1/3}\kappa^{-1}\|b\|^{2/3}_{L^{\infty}(0,T)} \right)T^{2}$, we can absorb the last two terms in \eqref{Source G}.\\
	Concerning the sixth and the second-to-last terms on the right-hand side of \eqref{trans2}, applying Young's inequality and $\hat{\xi}\geq \frac{C}{T^{2}}$, we obtain
	\begin{eqnarray}
		Cs \lambda^{2} \int_{0}^{T}\hat{\xi}(t)|\psi(L,t)||\psi_{x}(L,t)|&\leq& Cs \lambda^{3} \int_{0}^{T}\hat{\xi}(t)|\psi(L,t)|^{2} + \frac{s\lambda}{2} \int_{0}^{T}\hat{\xi}(t)|\psi_{x}(L,t)|^{2} \nonumber\\
		&\leq& Cs \lambda^{3}T^{4} \int_{0}^{T}\hat{\xi}^{3}(t)|\psi(L,t)|^{2}+ \frac{s\lambda}{2} \int_{0}^{T}\hat{\xi}(t)|\psi_{x}(L,t)|^{2} \nonumber \label{A}\\
	\end{eqnarray}
	and 
	\begin{eqnarray}
		&&C(\mu\kappa s\lambda+s\lambda^{2}) \int_{0}^{T} \hat{\xi}(t)|\psi(0,t)| |\psi_{x}(0,t)| \nonumber \\
		&&\leq  C(\mu^{2}\kappa^{2} s\lambda+s\lambda^{3}) \int_{0}^{T} \hat{\xi}(t)|\psi(0,t)|^{2}+ \frac{s\lambda}{2} \int_{0}^{T} \hat{\xi}(t)|\psi_{x}(0,t)|^{2} \nonumber\\
		&&\leq  C(\mu^{2}\kappa^{2} s\lambda+s\lambda^{3}) T^{4}\int_{0}^{T} \hat{\xi}^{3}(t)|\psi(0,t)|^{2}+ \frac{s\lambda}{2} \int_{0}^{T} \hat{\xi}(t)|\psi_{x}(0,t)|^{2}.  \label{B}
	\end{eqnarray}
	The first term on the right-hand side of \eqref{A} is absorbed if we take $s\geq s_{1}T^{2}$ and the second one is half a term on the left-hand side of \eqref{trans2}. 
	The same for the terms on right-hand side of \eqref{B} if we take $s\geq s_{1}(\mu\kappa+1)T^{2}$. The eleventh term on the right-hand side of \eqref{trans2} is absorbed  by the left-hand side of \eqref{trans2} by taking $s\geq (\mu\kappa)^{-1/2}T$.
	The last term on the right-hand side of \eqref{trans2} is absorbed  by the left-hand side of \eqref{trans2}. Other terms not included in the following are also absorbed. Consequently, we conclude 
	\begin{eqnarray}
		&&\|M_{1}\psi\|^{2}_{L^{2}(Q)}+\|M_{2}\psi\|^{2}_{L^{2}(Q)}+ \mu\kappa^{-1}\|N_{1}(\psi,\gamma)\|^{2}_{L^{2}(0,T)}+ \mu\kappa^{-1}\|N_{2}(\psi,\gamma)\|^{2}_{L^{2}(0,T)} \nonumber\\
		&&+ s^{3} \lambda^{4} \int_{Q} \xi^{3} |\psi|^{2} + s^{3}\lambda^{3}\int_{0}^{T}\hat{\xi}^{3}(t)\left(|\psi(0,t)|^{2}+ |\psi(L,t)|^{2}\right) \nonumber\\
		&&+s \lambda^{2}\int_{\Omega_{T}}\xi|\psi_{x}|^{2}  + s\lambda\int_{0}^{T}\hat{\xi}(t)\left(|\psi_{x}(0,t)|^{2}+ |\psi_{x}(L,t)|^{2}\right) \nonumber\\
		&&  \leq C\left( \int_{Q}e^{-2s\alpha}|f|^{2}+\mu\kappa^{-1}\int_{0}^{T}e^{-2s\hat{\alpha}}|g|^{2} +s^{3} \lambda^{4} \int_{\omega^{\prime}_T} \xi^{3} |\psi|^{2} + s \lambda^{2}\int_{\omega^{\prime}_{T}}\xi|\psi_{x}|^{2} \right),\nonumber\\
		\label{First conclusion}
	\end{eqnarray}
	for any $s\geq s_{1}\left(\left(1+(\mu\kappa)^{-1/2}\right)T+\left(1+\mu\kappa\right)T^{2}\right)$ and any $\lambda\geq\lambda_{1}$.\\
	\textbf{\underline{Step 4.} Indirect estimates and conclusion.} \\
	We start by adding integrals of $|\psi_{t}|^{2}$, $|\gamma_{t}|^{2}$ and $|\psi_{xx}|^{2}$ to the left-hand side of \eqref{First conclusion}, so that we can eliminate the last term in the right-hand side of \eqref{First conclusion}. Using \eqref{m1m2n1}, $\xi\geq \frac{C}{T^{2}}$, $s\geq s_{1}T^{2}$ and $|\alpha_{t}|\leqslant CT\xi^{2}$, we
	obtain
	%				$$M_{1}\psi:=-2s\lambda^{2}\psi\xi|\eta^{\prime}|^{2}-2s\lambda\xi\psi_{x}\eta^{\prime}+\psi_{t}$$
	%				$$-\psi_{t}=-2s\lambda^{2}\psi\xi|\eta^{\prime}|^{2}-2s\lambda\xi\psi_{x}\eta^{\prime}-M_{1}\psi$$
	\begin{eqnarray}
		s^{-1}\int_{Q}\xi^{-1}|\psi_{t}|^{2}dxdt
		&\leqslant&  C\left(s\lambda^{4}\int_{Q}\xi\psi^{2} +s \lambda^{2}\int_{Q}\xi\psi_{x}^{2}+ \|M_{1}\psi\|^{2}_{L^{2}(Q)}\right) \nonumber\\
		&\leqslant&  C\left(s^{3}\lambda^{4}\int_{Q}\xi^{3}\psi^{2} +s \lambda^{2}\int_{Q}\xi\psi_{x}^{2}+ \|M_{1}\psi\|^{2}_{L^{2}(Q)}\right) \label{m1}
	\end{eqnarray}
	%				for all $s\geq s_{1}T^{2}$ and all $\lambda\geq \lambda_{1}$.
	%				$$M_{2}\psi = s^{2}\lambda^{2}\psi\xi^{2}|\eta^{\prime}|^{2}+\psi_{xx} +s\psi\alpha_{t}$$
	%				$$-\psi_{xx} = s^{2}\lambda^{2}\psi\xi^{2}|\eta^{\prime}|^{2} +s\psi\alpha_{t}-M_{2}\psi$$
	and
	\begin{eqnarray}
		s^{-1}\int_{Q}\xi^{-1}|\psi_{xx}|^{2}
		&\leqslant&  C\left(s^{3}\lambda^{4}\int_{Q}\xi^{3}\psi^{2} +s T^{2}\int_{Q}\xi^{3}\psi^{2}+ \|M_{2}\psi\|^{2}_{L^{2}(Q)}\right) \nonumber\\
		&\leqslant&  C\left(s^{3}\lambda^{4}\int_{Q}\xi^{3}\psi^{2} + \|M_{2}\psi\|^{2}_{L^{2}(Q)}\right)\label{m2}
	\end{eqnarray}
	for all $s\geq s_{1}(T+T^{2})$. We also have\\
	%				$$N_{1}(\psi,\gamma) :=\gamma^{\prime}(t)-\kappa s\lambda\eta^{\prime}(0)\hat{\xi}\psi(0,t)$$
	%				$$\gamma^{\prime}(t)=N_{1}(\psi,\gamma)+\kappa s\lambda\eta^{\prime}(0)\hat{\xi}\psi(0,t)$$
	\begin{eqnarray}
		\mu\kappa^{-1}s^{-1}\int_{0}^{T}\hat{\xi}^{-1}|\gamma^{\prime}|^{2}
		&\leqslant&  C\left(\mu\kappa s\lambda^{2}\int_{0}^{T}\hat{\xi}\psi^{2}(0,t) + \mu\kappa^{-1}\|N_{1}(\psi,\gamma)\|^{2}_{L^{2}(0,T)}\right)\nonumber\\
		&\leqslant&  C\left(s^{3}\lambda^{3}\int_{0}^{T}\hat{\xi}^{3}\psi^{2}(0,t) +\mu\kappa^{-1}\|N_{1}(\psi,\gamma)\|^{2}_{L^{2}(0,T)}\right) \label{n1}
	\end{eqnarray}
	for all $s\geq s_{1}(\mu\kappa)^{1/2}T^{2}$.
	Consequently, we deduce from \eqref{First conclusion} and \eqref{m1}-\eqref{n1} that
	\begin{eqnarray}
		&&s^{-1}\int_{Q}\xi^{-1}|\psi_{t}|^{2} +s^{-1}\int_{Q}\xi^{-1}|\psi_{xx}|^{2}+ \mu\kappa^{-1}s^{-1}\int_{0}^{T}\hat{\xi}^{-1}|\gamma^{\prime}|^{2} \nonumber\\
		&&+ s^{3} \lambda^{4} \int_{Q} \xi^{3} |\psi|^{2} + s^{3}\lambda^{3}\int_{0}^{T}\hat{\xi}^{3}(t)\left(|\psi(0,t)|^{2}+ |\psi(L,t)|^{2}\right) \nonumber\\
		&&+s \lambda^{2}\int_{\Omega_{T}}\xi|\psi_{x}|^{2}  + s\lambda\int_{0}^{T}\hat{\xi}(t)\left(|\psi_{x}(0,t)|^{2}+ |\psi_{x}(1,t)|^{2}\right) \nonumber\\
		&&  \leq C\left( \int_{Q}e^{-2s\alpha}|f|^{2}+\mu\kappa^{-1}\int_{0}^{T}e^{-2s\hat{\alpha}}|g|^{2} +s^{3} \lambda^{4} \int_{\omega^{\prime}_T} \xi^{3} \psi^{2} + s \lambda^{2}\int_{\omega^{\prime}_{T}}\xi|\psi_{x}|^{2} \right)
		\label{Second conclusion}
	\end{eqnarray}
	for any $s\geq s_{1}\left((1+(\mu\kappa)^{-1/2})T+\left(1+\mu\kappa+(\mu\kappa)^{1/2}\right)T^{2}\right)$ and any $\lambda\geq\lambda_{1}$.\\
	\par 
	As usual, to eliminate the last term on the right-hand side of \eqref{Second conclusion}, let us introduce $\theta\in \mathcal{C}^{2}(\omega)$ a positive cut-off function such that $\theta=1$ in $\omega^{\prime}$, an integration by parts and the Cauchy-Schwarz inequality as in \cite{guerrero2007singular}, we obtain
	\begin{eqnarray*}
		Cs\lambda^{2}\int_{\omega^{\prime}_{T}}|\psi_{x}|^{2}\xi  &=& Cs\lambda^{2}\int_{\omega^{\prime}_{T}}\theta|\psi_{x}|^{2}\xi \leq Cs\lambda^{2}\int_{\omega_{T}}\theta|\psi_{x}|^{2}\xi \\
		&\leq & \frac{1}{2}\left(s^{-1}\int_{Q}\xi^{-1}|\psi_{xx}|^{2}+s\lambda^{2}\int_{Q}\xi|\psi_{x}|^{2}\right) + Cs^{3}\lambda^{4}\int_{\omega_{T}}\xi^{3}|\psi|^{2}.
	\end{eqnarray*}
	Combining this last estimate with \eqref{Second conclusion}, we deduce  
	\begin{eqnarray}
		&&s^{-1}\int_{Q}\xi^{-1}|\psi_{t}|^{2} +s^{-1}\int_{Q}\xi^{-1}|\psi_{xx}|^{2}+ \mu\kappa^{-1}s^{-1}\int_{0}^{T}\hat{\xi}^{-1}|\gamma^{\prime}|^{2} \nonumber\\
		&&+ s^{3} \lambda^{4} \int_{Q} \xi^{3} |\psi|^{2}  + s^{3}\lambda^{3}\int_{0}^{T}\hat{\xi}^{3}\left(|\psi(0,t)|^{2}+ |\psi(L,t)|^{2}\right) \nonumber\\
		&&+s \lambda^{2}\int_{\Omega_{T}}\xi|\psi_{x}|^{2}  + s\lambda\int_{0}^{T}\hat{\xi}\left(|\psi_{x}(0,t)|^{2}+ |\psi_{x}(L,t)|^{2}\right) \nonumber\\
		&&  \leq C\left( \int_{Q}e^{-2s\alpha}|f|^{2}+\mu\kappa^{-1}\int_{0}^{T}e^{-2s\hat{\alpha}}|g|^{2} +s^{3} \lambda^{4} \int_{\omega_T} \xi^{3} |\psi|^{2}   \right)\label{last estimate}
	\end{eqnarray}
	for any $s\geq s_{1}\left((1+(\mu\kappa)^{-1/2})T+\left(1+\mu\kappa+(\mu\kappa)^{1/2}\right)T^{2}\right)$ and any $\lambda\geq\lambda_{1}$.
	\par 
	Finally, we return to our original functions, which were as follows $(\varphi, \rho)=(e^{s\alpha}\psi, e^{s\hat{\alpha}}\gamma)$ and using the elementary identities \eqref{elementary identities}. For $\varphi_{x}$, using that
	\begin{eqnarray*}
		\varphi_{x}=e^{s\alpha}(\psi_{x}-s\lambda\eta^{\prime}\xi\psi),
	\end{eqnarray*}
	we find
	\begin{eqnarray}
		s\lambda^{2}\int_{Q}e^{-2s\alpha}\xi|\varphi_{x}|^{2}\leq s\lambda^{2}\int_{Q}\xi|\psi_{x}|^{2}+s^{3}\lambda^{4}\int_{Q}\xi^{3}|\psi|^{2} \label{reverse 1}
	\end{eqnarray}
	and 
	\begin{eqnarray}
		s\lambda\int_{0}^{T}e^{-2s\hat{\alpha}}\hat{\xi}|\varphi_{x}(0,t)|^{2}\leq s\lambda\int_{0}^{T}\hat{\xi}|\psi_{x}(0,t)|^{2}+s^{3}\lambda^{3}\int_{Q}\hat{\xi}^{3}|\psi(0,t)|^{2}. \label{reverse 2}
	\end{eqnarray}
	For $\varphi_{xx}$, using  the identity
	\begin{eqnarray*}
		\psi_{xx}= e^{-s \alpha} \varphi_{xx}+2s\lambda\xi\psi_{x}\eta^{\prime}-s^{2}\lambda^{2}\xi^{2}\psi|\eta^{\prime}|^{2}+s\lambda^{2}\xi\psi|\eta^{\prime}|^{2}+s\lambda\xi\psi\eta^{\prime\prime},
	\end{eqnarray*}
	we obtain
	\begin{eqnarray}
		&&s^{-1}\int_{Q}e^{-2s\alpha}\xi^{-1}|\varphi_{xx}|^{2}\leq s^{-1}\int_{Q}\xi^{-1}|\psi_{xx}|^{2}dxdt+s\lambda^{2}\int_{Q}\xi|\psi_{x}|^{2}  \nonumber\\
		&& + s^{3}\lambda^{4}\int_{Q}\xi^{3}|\psi|^{2}+ s\lambda^{4}\int_{Q}\xi|\psi|^{2}+ s\lambda^{2}\int_{Q}\xi|\psi|^{2} \nonumber\\
		&&\leq C\left(s^{-1}\int_{Q}\xi^{-1}|\psi_{xx}|^{2}+s\lambda^{2}\int_{Q}\xi|\psi_{x}|^{2} + s^{3}\lambda^{4}\int_{Q}\xi^{3}|\psi|^{2}\right)  \label{reverse 3}
	\end{eqnarray}
	for any $s\geq s_{1}T^{2}$.
	Finally, for $\varphi_{t}$ and $\rho_{t}$, we obtain
	\begin{eqnarray}
		s^{-1}\int_{Q}e^{-2s\alpha}\xi^{-1}|\varphi_{t}|^{2} &\leq& s^{-1}\int_{Q}\xi^{-1}|\psi_{t}|^{2}+CsT^{2}\int_{Q}\xi^{3}|\psi|^{2} \nonumber\\
		&\leq& s^{-1}\int_{Q}\xi^{-1}|\psi_{t}|^{2}+Cs^{3}\lambda^{4}\int_{Q}\xi^{3}|\psi|^{2} \label{reverse 4}
	\end{eqnarray}
	for any $s\geq s_{1}T$ and 
	\begin{eqnarray}
		\mu\kappa^{-1}s^{-1}\int_{0}^{T}e^{-2s\hat{\alpha}}\hat{\xi}^{-1}|\rho^{\prime}|^{2} &\leq& \mu\kappa^{-1}s^{-1}\int_{0}^{T}\hat{\xi}^{-1}|\gamma^{\prime}|^{2}+C\mu\kappa^{-1}sT^{2}\int_{0}^{T}\hat{\xi}^{3}|\gamma|^{2} \nonumber\\
		&\leq& s^{-1}\int_{0}^{T}\hat{\xi}^{-1}|\gamma^{\prime}|^{2}+Cs^{3}\lambda^{3}\int_{0}^{T}\hat{\xi}^{3}|\psi(0,t)|^{2} \label{reverse 5}
	\end{eqnarray}
	for any $s\geq s_{1}(\mu\kappa)^{-1}T$. Consequently, from \eqref{last estimate} and \eqref{reverse 1}-\eqref{reverse 5}, we can derive \eqref{Carleman estimate}.
\end{proof}
\par 
Using Carleman estimate \eqref{Carleman estimate}, the following observability inequality holds.
\begin{corollary}
	There is a constant $C_{ob}>0$ such that for all $(\varphi_{T},\rho_{T})\in \mathcal{H}$ the weak solution $(
	\varphi, \rho)$ of the adjoint system \eqref{s2} satisfies
	\begin{eqnarray}
		\int_{0}^{L}|\varphi(x,0)|^{2}\d x+ \frac{\mu}{\kappa}|\rho(0)|^{2}\leq C_{ob}\int_{\omega_T}|\varphi|^{2}\d x\d t. \label{observability inequality1}
	\end{eqnarray}
\end{corollary}
\begin{proof}
	By density we can thus restrict ourselves to final values $(\varphi_{T},\rho_{T})\in \mathcal{H}^{1}$, so that $(
	\varphi, \rho)$ is a strong solution of the adjoint system \eqref{s2}. Applying Proposition \ref{P4} to system \eqref{s2} with $f=-a(x,t)\varphi$ and $g=-\mu^{-1}c(t)\varphi(0,t)-\mu^{-1}\kappa\displaystyle\int_{0}^{L}b(x,t)\varphi(x,t)\d x$ leads to 
	\begin{eqnarray}
		&&s^{3}\lambda^{4}\int_{Q}e^{-2s\alpha}\xi^{3}|\varphi|^{2}\d x\d t +s^{3}\lambda^{3}\int_{0}^{T}e^{-2s\hat{\alpha}}\hat{\xi}^{3}|\varphi(0,t)|^{2}\d t  \nonumber\\
		&& \leq C\left(s^{3}\lambda^{4}\int_{\omega_{T}}e^{-2s\alpha}\xi^{3}|\varphi|^{2}\d x\d t + \int_{Q}e^{-2s\alpha}|a(x,t)\varphi|^{2}\d x\d t \right. \nonumber\\
		&& \left. + \mu^{-1}\kappa^{-1}\int_{0}^{T}e^{-2s\hat{\alpha}}\left|c(t)\varphi(0,t)\right|^{2}\d t +\mu^{-1}\kappa\int_{0}^{T}e^{-2s\hat{\alpha}}\left|\int_{0}^{L}b(x,t)\varphi(x,t)dx\right|^{2}\d t  \right). \nonumber
	\end{eqnarray}
	Using the Cauchy Schwarz inequality in space for the nonlocal term and $\xi\geq \frac{C}{T^{2}}$, we can absorb the last three terms of the right-hand side of the above inequality by taking $s\geq s_1\left( \|a\|_{L^{\infty}(Q)}^{2/3}+ (\mu\kappa)^{-1/3}\|c\|_{L^{\infty}(0,T)}^{2/3}+ \mu^{-1/3}\kappa^{1/3}\|b\|_{L^{\infty}(Q)}^{2/3}\right)T^{2}$. As a result, we obtain the following.
	\begin{eqnarray*}
		\int_{Q}e^{-2s\alpha}\xi^{3}|\varphi|^{2}\d x\d t +\frac{\mu^{-1}}{\kappa}\int_{0}^{T}e^{-2s\hat{\alpha}}\hat{\xi}^{3}|\varphi(0,t)|^{2}\d t \leq C\int_{\omega_{T}}e^{-2s\alpha}\xi^{3}|\varphi|^{2}\d x\d t
	\end{eqnarray*}
	for any $s\geq s_{2}(T+T^{2})$ and any $\lambda\geq \lambda_{2}$, where $C, s_2$ and $\lambda_{2}$ also depend on $\mu, \kappa, a, b$ and $c$. Let us fix $s= s_{2}(T+T^{2})$ and $\lambda=\lambda_{2}$, we get 
	\begin{eqnarray}
		\int_{(\frac{T}{4},\frac{3T}{4})\times (0,L)}|\varphi|^{2}\d x\d t +\frac{\mu}{\kappa}\int_{\frac{T}{4}}^{\frac{3T}{4}}|\rho(t)|^{2}\d t \leq C\int_{\omega_{T}}|\varphi|^{2}\d x\d t. \label{es}
	\end{eqnarray}
	As usual, the estimate \eqref{es} and the dissipation properties of the solutions of \eqref{s2}, will lead to \eqref{observability inequality1}.
\end{proof}
%			\begin{proof}
	%				A classical and usual argument allows us to obtain this obervability inequality from Carleman's estiamtion \eqref{Carleman estimate} as in \cite{bibid}.
	%			\end{proof}
%			\subsection{Proof of Theorem \ref{Main 2}} In this section we apply the Carleman estimate to show null controllability of \eqref{s0} and its generalizations.  *************
%			By a standard duality argument (Hilbert's uniqueness method), the null controllability of \eqref{s0} is equivalent to the observability inequality \eqref{observability inequality} of the adjoint system \eqref{s1}. \\
%			The observability inequality \eqref{observability inequality1} implies the null controllability
%			of \eqref{s0}.
\subsection{Proof of Theorem \ref{Main 2}}
For abstract linear control systems, it is well known that their controllability is equivalent to the observability inequality of their adjoint system. In addition, if control exists, it is certainly not unique. This is why, to prove controllability results, it is useful to specify such a control. For instance, the HUM approach consists in finding the control with
the minimal $L^{2}(\omega_{T})$-norm as specified in the Proof of Theorem \ref{Main 2}.
\begin{notation}
	Denote systems \eqref{s1} and \eqref{s2} by $\mathcal{S}(v,y_{0},z_{0})$ and $\mathcal{S}^{*}(\varphi_{T},\rho_{T})$ respectively.
\end{notation}
\begin{proof}[Proof of Theorem~{\upshape\ref{Main 2}}]
	We split the proof into two steps. First, we construct a sequence of controls $v^{\varepsilon}\in L^{2}(\omega_{T})$ with $\varepsilon>0$ that yield the approximate null controllability of \eqref{s0}. Second, we conclude by passing to the limit when $\varepsilon$ tends to zero.\\
	\textbf{Step 1:} Fix $(y_{0},z_{0})\in \mathcal{H}$ an initial state to be controlled and let us consider, for all $\varepsilon>0$, the following optimal control problem
	\begin{eqnarray*}
		\min_{v\in L^{2}(\omega_{T})}\mathbf{F}_{\varepsilon}(v),
	\end{eqnarray*}
	where
	\begin{eqnarray*}
		\mathbf{F}_{\varepsilon}(v)=\frac{1}{2}\int_{\omega_{T}}|v|^{2}dxdt +\frac{1}{2\varepsilon}\|(y(\cdot,T), z(T))\|^{2}_{\mathcal{H}}
	\end{eqnarray*}
	and $(y,z)$ is the solution of $\mathcal{S}(v,y_{0},z_{0})$.
	For any $\varepsilon>0$, the functional $\mathbf{F}_{\varepsilon}$ has a unique minimizer $v_{\varepsilon}\in L^{2}(\omega_{T})$. This is
	due to the fact that $\mathbf{F}_{\varepsilon}$ is strictly convex, continuous and coercive. This minimizer is characterized by the following optimality condition of first order (Euler-Lagrange equation)
	%\begin{eqnarray}
	%	\int_{\omega_{T}}v_{\varepsilon}v dxdt +\frac{1}{\varepsilon}\left\langle \mathcal{S}_{T}(v_{\varepsilon},y_{0},z_{0}), \mathcal{S}_{T}(v,0,0) \right\rangle_{\mathcal{H}}=0\quad \forall v\in L^{2}(\omega_{T}). \label{E-L}
	%\end{eqnarray}
	\begin{eqnarray}
		\int_{\omega_{T}}v_{\varepsilon}v dxdt +\frac{1}{\varepsilon}\left\langle (y^{\varepsilon}(\cdot,T), z^{\varepsilon}(T)), (y^{0}(\cdot,T), z^{0}(T)) \right\rangle_{\mathcal{H}}=0\quad \forall v\in L^{2}(\omega_{T}), \label{E-L}
	\end{eqnarray}
	where $(y^{\varepsilon},z^{\varepsilon})$ and $(y^{0},z^{0})$ are respectively the solutions of $\mathcal{S}(v^{\varepsilon},y_{0},z_{0})$ and $\mathcal{S}(v,0,0)$.\\
	On the other hand, we have
	\begin{eqnarray*}
		\int_{\omega_{T}}v\varphi dxdt=\left\langle (y^{0}(\cdot,T), z^{0}(T)), (\varphi_{T},\rho_{T}) \right\rangle_{\mathcal{H}} ,\quad \forall (\varphi_{T},\rho_{T})\in \mathcal{H},
	\end{eqnarray*}
	where $(\varphi,\rho)$ is the solution of $\mathcal{S}^{*}(\varphi_{T},\rho_{T})$.
	In particular, choosing $(\varphi_{T},\rho_{T})=\frac{1}{\varepsilon}(y^{\varepsilon}(\cdot,T),z^{\varepsilon}(T))$, we obtain 
	\begin{eqnarray*}
		\int_{\omega_{T}}v\varphi^{\varepsilon} dxdt=\frac{1}{\varepsilon}\left\langle (y^{\varepsilon}(\cdot,T), z^{\varepsilon}(T)), (y^{0}(\cdot,T), z^{0}(T)) \right\rangle_{\mathcal{H}},  \quad \forall v\in L^{2}(\omega_{T}),
	\end{eqnarray*}
	where $(\varphi^{\varepsilon},\rho^{\varepsilon})$ is the solution of $\mathcal{S}^{*}\left(\frac{1}{\varepsilon}y^{\varepsilon}(\cdot,T), \frac{1}{\varepsilon}z^{\varepsilon}(T) \right)$. This, combined with \eqref{E-L}, provides 
	\begin{eqnarray}
		v_{\varepsilon}=-\mathds{1}_{\omega}\varphi^{\varepsilon}. \label{Control}
	\end{eqnarray}
	Using \eqref{Control} and applying the duality relation between $\mathcal{S}^{*}\left(\frac{1}{\varepsilon}y^{\varepsilon}(\cdot,T), \frac{1}{\varepsilon}z^{\varepsilon}(T) \right)$ and $\mathcal{S}(v^{\varepsilon},y_{0},z_{0})$, we obtain
	%\begin{eqnarray}
	%	\int_{\omega_{T}}v\varphi dxdt=\left\langle \begin{pmatrix}
		%		y(\cdot,T)  \\ z(T)
		%	\end{pmatrix}, \begin{pmatrix}
		%		\varphi_{T}  \\ \rho_{T}
		%	\end{pmatrix} \right\rangle_{\mathcal{H}}-\left\langle \begin{pmatrix}
		%		y_{0}  \\ z_{0}
		%	\end{pmatrix}, \begin{pmatrix}
		%		\varphi(\cdot,0)  \\ \rho(0)
		%	\end{pmatrix} \right\rangle_{\mathcal{H}}. 
	%\end{eqnarray}
	\begin{eqnarray*}
		-\int_{\omega_{T}}|v^{\varepsilon}|^{2} dxdt=\frac{1}{\varepsilon}\|(y^{\varepsilon}(\cdot,T), z^{\varepsilon}(T)\|^{2}_{\mathcal{H}}-\left\langle (y_{0}, z_{0})
		, (\varphi^{\varepsilon}(\cdot,0), \rho^{\varepsilon}(0))\right\rangle_{\mathcal{H}}. 
	\end{eqnarray*}
	Thanks to Cauchy-Schwarz inequality and the observability inequality \eqref{observability inequality1},
	\begin{eqnarray*}
		\frac{1}{\varepsilon}\|(y^{\varepsilon}(\cdot,T), z^{\varepsilon}(T)\|^{2}_{\mathcal{H}}+\int_{\omega_{T}}|v^{\varepsilon}|^{2} dxdt
		%&=&\left\langle (y_{0}, z_{0}), (\varphi^{\varepsilon}(\cdot,0), \rho^{\varepsilon}(0))\right\rangle_{\mathcal{H}}\\
		\leq C_{op}^{1/2}\|(y_{0}, z_{0}\|_{\mathcal{H}}\left(\int_{\omega_{T}}|v^{\varepsilon}|^{2} dxdt\right)^{1/2}.
	\end{eqnarray*}
	Consequently
	\begin{eqnarray}
		\int_{\omega_{T}}|v^{\varepsilon}|^{2} dxdt\leq C_{op}\|(y_{0}, z_{0}\|^{2}_{\mathcal{H}} \label{control-state A}
	\end{eqnarray}
	and 
	\begin{eqnarray}
		\|(y^{\varepsilon}(\cdot,T), z^{\varepsilon}(T)\|^{2}_{\mathcal{H}}
		\leq \varepsilon C_{op}\|(y_{0}, z_{0}\|^{2}_{\mathcal{H}}. \label{bounded}
	\end{eqnarray}
	\textbf{Step 2:} Since $(v^{\varepsilon})$ is bounded in $L^{2}(\omega_{T})$, it possesses a (weakly) convergent subsequence to certain $v\in L^{2}(\omega_{T})$. Using classical
	parabolic estimates we deduce that, at least for a subsequence,
	\begin{eqnarray*}
		(y^{\varepsilon},z^{\varepsilon})\longrightarrow (y,z)\;\mbox{weakly in}\; L^{2}(0,T; \mathcal{H}^{1})\cap H^{1}(0,T; (\mathcal{H}^{1})^{\prime}),
	\end{eqnarray*}
	where $(y,z)$ is the weak solution of $\mathcal{S}(v,y_{0},z_{0})$. In particular, this gives a weak convergence of $(y^{\varepsilon}(\cdot,T), z^{\varepsilon}(T))$ to $(y(\cdot,T), z(T))$ in $\mathcal{H}$, which combined with \eqref{bounded}, yields $(y(\cdot,T),z(T))=(0,0)$. Moreover, \eqref{control-state} is a consequence of \eqref{control-state A}.
\end{proof}
\subsection{Proof of Theorem \ref{Main 1}}
In this subsection, we give the proof of Theorem \ref{Main 1}.
\begin{proof}[Proof of Theorem~{\upshape\ref{Main 1}}] Let us consider a nonempty open set $\omega\subset (\ell, L)$ for $L>\ell$ and $\tilde{y_{0}}\in H^{1}(0,L)$, such that $\tilde{y_{0}}=y_{0}$ in $(0,\ell)$. 
	Now, we consider the following interior control problem
	\begin{equation}  \label{BC2}
		\left\{
		\begin{aligned}
			&y_{t}-y_{xx} + a(x,t)y+b(x,t)z(t)=\mathds{1}_{\omega}v
			& & \text {in}\; Q_{L}, \\
			& z^{\prime}(t)+c(t)z(t)-\kappa y_{x}(0,t)
			=0 & & \text {in}\;(0,T), \\
			& y(0,t)=\mu z(t)  & & \text {in}\;(0,T), \\
			& y_{x}(L,t) =0 & & \text {in}\;(0,T), \\
			& y(\cdot,0)=\tilde{y_{0}}  & & \text {in}\;(0,L), \\
			&z(0)=z_{0},   \\	
		\end{aligned}
		\right.
	\end{equation}
	The interior null controllability
	result in Theorem \ref{Main 2} gives a control $v\in L^{2}(\omega_{T})$ such that the solution $(y,z)\in \mathcal{E}(0,L)$ of \eqref{BC2} satisfies 
	\begin{eqnarray*}
		y(\cdot,T)=0 \;\text{in}\; (0,L) \quad \mbox{and} \quad z(T)=0. 
	\end{eqnarray*}
	We then define the control $u=y_{x}(\ell,\cdot)\in H^{1/4}(0,T)$ due to $y_{|_{(0,\ell)}}\in L^{2}(0,T;H^{2}(0,\ell))\cap H^{1}(0,T; L^{2}(0,\ell))$ and \cite[Theorem 2.1]{lions1972non}. Consequently $(y_{|_{(0,\ell)}},z)\in \mathcal{E}(0,\ell)$ is a solution of \eqref{s0} associated with control $u$ and satisfies \eqref{Null}.
\end{proof}
\section{Some numerical results and experiments} \label{sec4}
\subsection{Algorithm for calculating HUM controls}
In this section, we look at a numerical algorithm to calculate HUM controls. This method uses a penalized HUM approach and a conjugate gradient algorithm (CG for short). We refer to \cite{glowinskiexact} and \cite{boyer2013penalised} for more details on such a method.
As the proof of Proposition \ref{Main 2} shows, the solution associated with the minimizer of $\textbf{F}_{\varepsilon}$ will be an approximation of the null controllability problem. Applying general results of the Fenchel-Rockafellar theory, we can construct an associated dual
problem as follows.
\begin{eqnarray*}
	\mathbf{J}_{\varepsilon}(\varphi_{T},\rho_{T})=\frac{1}{2}\int_{\omega_{T}}|\varphi|^{2}\d x\d t +\frac{\varepsilon}{2}\|(\varphi_{T},\rho_{T})\|^{2}_{\mathcal{H}}+\langle (y_{0},z_{0}),(\varphi(\cdot,0),\rho(0))\rangle_{\mathcal{H}},
\end{eqnarray*}
where $(\varphi,\rho)$ is the solution of $\mathcal{S}^{*}(\varphi_{T},\rho_{T})$.
The duality properties between these two functional are consequence of the general results mentioned above. More precisely, we have the following result, see \cite[Proposition 1.5]{boyer2013penalised}.
\begin{proposition}
	For any $\varepsilon>0$, the minimizers $v^{\varepsilon}$ and $(\varphi^{\varepsilon}_{T},\rho^{\varepsilon}_{T})$ of the functionals $\mathbf{F}_{\varepsilon}$ and $\mathbf{J}_{\varepsilon}$ respectively, are
	related through the formulas
	\begin{eqnarray*}
		v_{\varepsilon}(\cdot,t)=\mathds{1}_{\omega}\varphi^{\varepsilon}(\cdot,t)\quad\mbox{for a.e.}\; t\in (0,T),
	\end{eqnarray*}
	where $(\varphi^{\varepsilon},\rho^{\varepsilon})$ is the solution of $\mathcal{S}^{*}(\varphi^{\varepsilon}_{T},\rho^{\varepsilon}_{T})$. Moreover, we have
	\begin{eqnarray*}
		\min_{v\in L^{2}(\omega_{T})}\mathbf{F}_{\varepsilon}(v)=\mathbf{F}_{\varepsilon}(v^{\varepsilon})=-\mathbf{J}_{\varepsilon}(\varphi^{\varepsilon}_{T},\rho^{\varepsilon}_{T})=-\min_{(\varphi_{T},\rho_{T})\in\mathcal{H}}\mathbf{J}_{\varepsilon}(\varphi_{T},\rho_{T}),
	\end{eqnarray*}
\end{proposition}
In the following, we will proceed with the dual problem.
The minimizer of $\textbf{J}_{\varepsilon}$ is characterized by the Euler-Lagrange equation
\begin{eqnarray}
	\int_{\omega_{T}}\varphi^{\varepsilon}\varphi\d x\d t +\varepsilon \langle (\varphi_{T}^{\varepsilon},\rho_{T}^{\varepsilon}),(\varphi_{T},\rho_{T})\rangle_{\mathcal{H}} +\langle (y_{0},z_{0})(\varphi(\cdot,0),\rho(0))\rangle_{\mathcal{H}}=0 \label{EL}
\end{eqnarray}
for all $(\varphi_{T},\rho_{T})\in\mathcal{H}$, where $(\varphi,\rho)$ is the solution of $\mathcal{S}^{*}(\varphi_{T},\rho_{T})$. 

Let us now define the linear operator $\Lambda$, usually referred to as the Gramian operator, as follows
\begin{eqnarray*}
	\Lambda(\varphi_{T},\rho_{T})=(y(\cdot,T),z(T)), 
\end{eqnarray*}
where $(y,z)$ is the solution of $\mathcal{S}(\varphi\mathds{1}_{\omega},0,0)$, while $(\varphi,\rho)$ is the solution of $\mathcal{S}^{*}(\varphi_{T},\rho_{T})$. The duality relation \eqref{dr} applied to $\mathcal{S}(\varphi^{\varepsilon}\mathds{1}_{\omega},0,0)$ and $\mathcal{S}^{*}(\varphi_{T},\rho_{T})$, yields 
\begin{eqnarray}
	\int_{\omega_{T}}\varphi\varphi^{\varepsilon}\d x\d t=\langle \Lambda(\varphi^{\varepsilon}_{T},\rho^{\varepsilon}_{T}),(\varphi_{T},\rho_{T})\rangle_{\mathcal{H}}. \label{EL1}
\end{eqnarray}
Again applying the duality relation between $\mathcal{S}(0,y_{0},z_{0})$ and $\mathcal{S}^{*}(\varphi_{T},\rho_{T})$, we obtain
\begin{eqnarray}
	\langle (y_{0},z_{0}),(\varphi(\cdot,0),\rho(0))\rangle_{\mathcal{H}}=\left\langle (
	y(\cdot,T), z(T))
	, (\varphi_{T},\rho_{T}) \right\rangle_{\mathcal{H}}, \label{EL2}
\end{eqnarray}
where $(y,z)$ is the solution of $\mathcal{S}(0,y_{0},z_{0})$.
By injecting \eqref{EL1} and \eqref{EL2} into \eqref{EL},  we obtain the following linear equation:
\begin{eqnarray*}
	(\Lambda+\varepsilon I_{\mathcal{H}})(\varphi^{\varepsilon}_{T},\rho^{\varepsilon}_{T})=-(
	y(\cdot,T), z(T)).
\end{eqnarray*}
To resolve this operator equation, we propose the following CG algorithm.
\begin{algorithm}[H]
	\caption{HUM combined with CG Algorithm}\label{algo1}
	\begin{algorithmic}[1]
		\State Set $k=0$ and choose an initial guess $f^{0}\in\mathcal{H}$
		\State Compute $(\varphi^{0}, \rho^{0})$ the solution of $\mathcal{S}^{*}(f^{0})$
		\State Set $u^{0}=\varphi^{0}\mathds{1}_{\omega}$ 
		\State Compute $(y^0, z^0)$ the solution of $\mathcal{S}(u^{0},y_{0},z_{0})$ 
		\State Compute $g^{0}=\varepsilon f^{0}+(y^{0}(\cdot,T),z^{0}(T))$
		\If{$\frac{\|g^{0}\|_{\mathcal{H}}}{\|f^{0}\|_{\mathcal{H}}}\leq \mbox{tol}$}\label{algln2}
		\State $f=f^{0}$
		\State $u=u^{0}$
		\Else
		\State $w^{0}=g^{0}$
		\EndIf
		\State For $k=1,2,\cdots,$ assuming that $f^{k}, g^{k}$ and $w^{k}$ are known, compute $f^{k+1}, g^{k+1}$ and $w^{k+1}$ as follows:
		\State \label{line}Compute $(\varphi^{k},\rho^{k})$ the solution of $\mathcal{S}^{*}(w^{k})$
		\State Set $u^{k}=\varphi^{k}\mathds{1}_{\omega}$ 
		\State Compute $(y^{k},z^{k})$ the solution of $S(u^{k},0,0)$
		\State Compute \begin{eqnarray*}
			&&\bar{g}^{k}= \varepsilon w^{k}+(y^{k}(\cdot,T),z^{k}(T))\\
			&&\varrho^{k}= \frac{\|g^{k}\|^{2}_{\mathcal{H}}}{\langle \bar{g}^{k}, w^{k}\rangle}_\mathcal{H}\\
			&&f^{k+1}= f^{k}-\varrho^{k}w^{k}\\
			&&g^{k+1}= g^{k}-\varrho^{k}\bar{g}^{k}
		\end{eqnarray*}
		\If{$\frac{\|g^{k+1}\|_{\mathcal{H}}}{\|g^{0}\|_{\mathcal{H}}}\leq \mbox{tol}$}
		\State Determine $(\varphi,\rho)$ the solution of $\mathcal{S}^{*}(f^{k+1})$
		\State Set $u=\varphi\mathds{1}_{\omega}$ as the minimizer of $\textbf{F}_{\varepsilon}$
		\Else 
		\State Compute
		\begin{eqnarray*}
			&&\gamma^{k}=\frac{\|g^{k+1}\|^{2}_{\mathcal{H}}}{\|g^{k}\|^{2}_{\mathcal{H}}}\\
			&&w^{k+1}= g^{k+1}+\gamma^{k}w^{k}\\
			&&k \gets k + 1\\
			&& \mbox{\Return to line \ref{line}}
		\end{eqnarray*}
		\EndIf
	\end{algorithmic}
\end{algorithm}
\subsection{Some numerical experiments}
We will solve numerically the null controllability problem
for \eqref{s0} and \eqref{s1}, and we will check that the previous CG Algorithm converge satisfactorily in several particular cases.
\subsubsection{Test 1}  
The CG Algorithm has been applied to the solution of the null controllability
problem for \eqref{s1} with the following data:
\begin{itemize}
	\item $L=1$, $\omega=(0.3,0.7)$, $T=0.6$.
	\item $y_{0}(x)=-10\sin(\pi x)$, $z_{0}=0$.
	\item $a=1, b=0, c=1$ and $\mu=\kappa=1$.
\end{itemize}
\par
For our computations, we take $Nx=30$ for the spatial mesh parameter. The initial guess in the algorithm is taken as $f^0=(0.4\sin(\pi x),0)$. We also choose $\varepsilon= 10^{-4}$ and the stopping parameter tol$=10^{-3}$ for the plots.
\begin{figure}[H]
	\begin{minipage}[c]{.3\linewidth}
		\begin{center}
			\includegraphics[width=4cm]{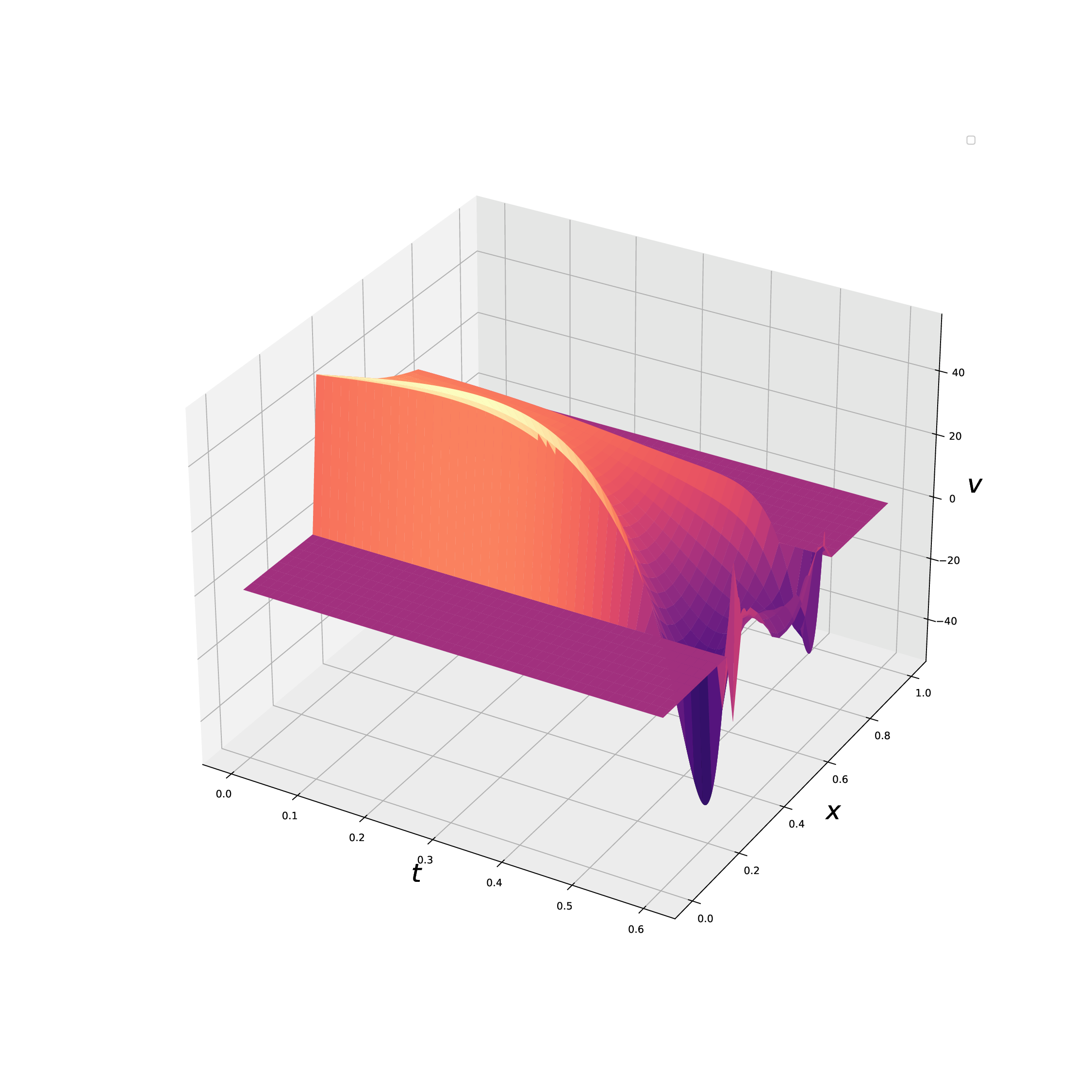}
		\end{center}
		\title{(a) Computed control.}
	\end{minipage} \hfill
	\begin{minipage}[c]{.3\linewidth}
		\begin{center}
			\includegraphics[width=4cm]{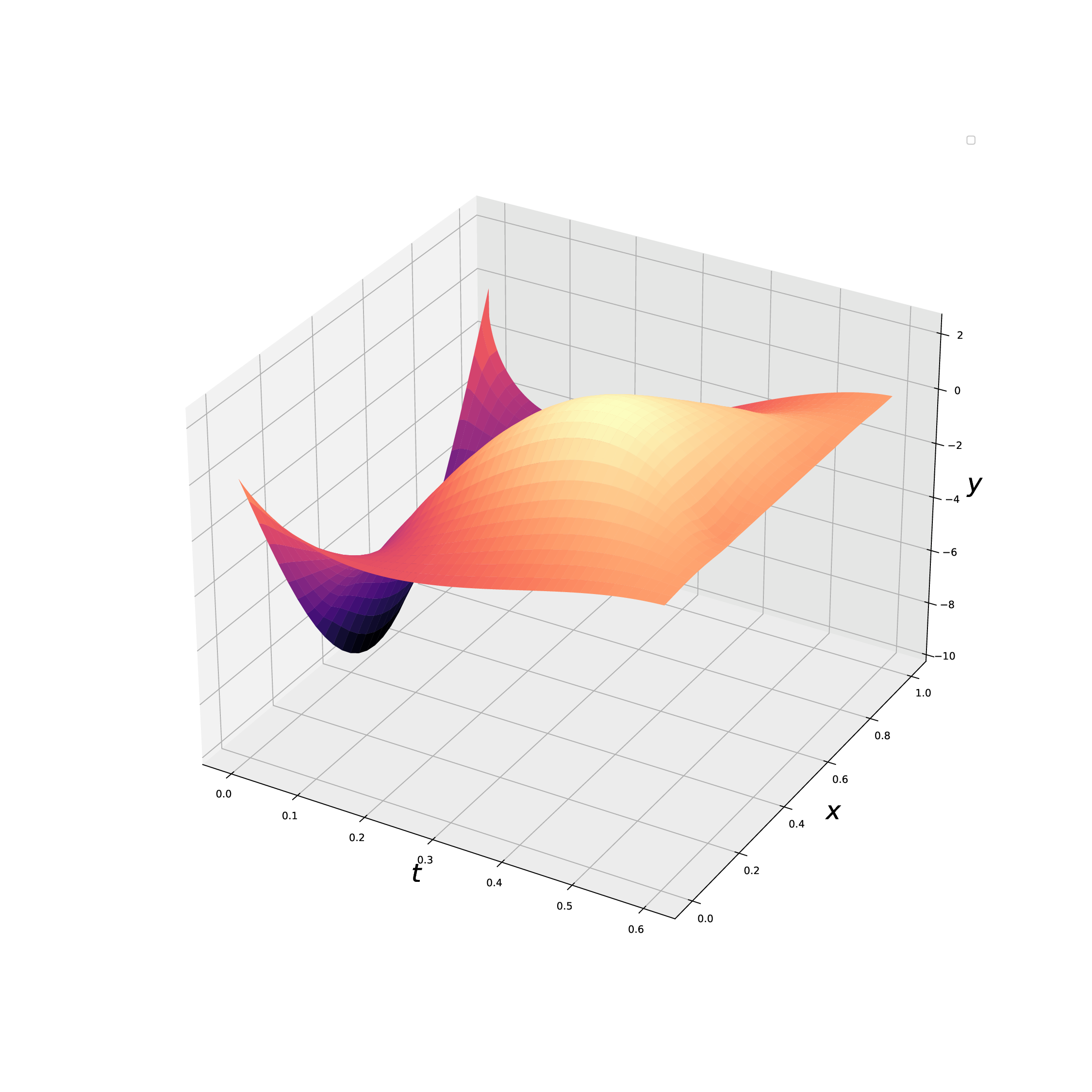}
		\end{center}
		\title{(b) Associated state $y$.}
	\end{minipage}
	\hfill
	\begin{minipage}[c]{.3\linewidth}
		\begin{center}
			\includegraphics[width=4cm]{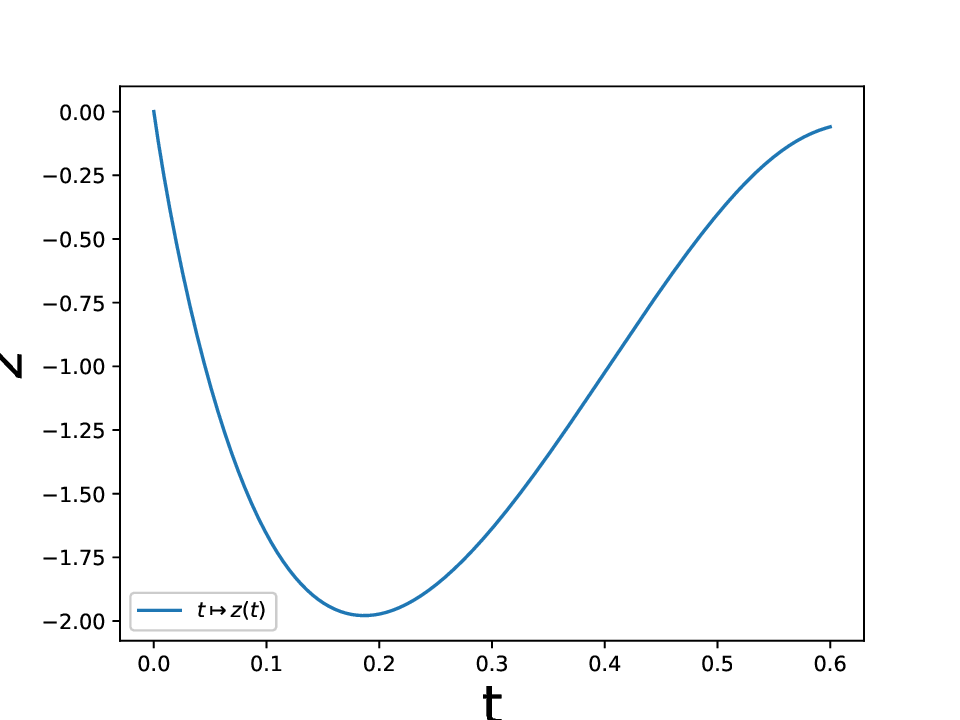}
		\end{center}
		\title{(c) Associated state $z$.}
	\end{minipage}
	\caption{The computed control and the associated state.}
\end{figure}

\begin{figure}[H]
	\begin{minipage}[c]{.4\linewidth}
		\begin{center}
			\includegraphics[width=5cm]{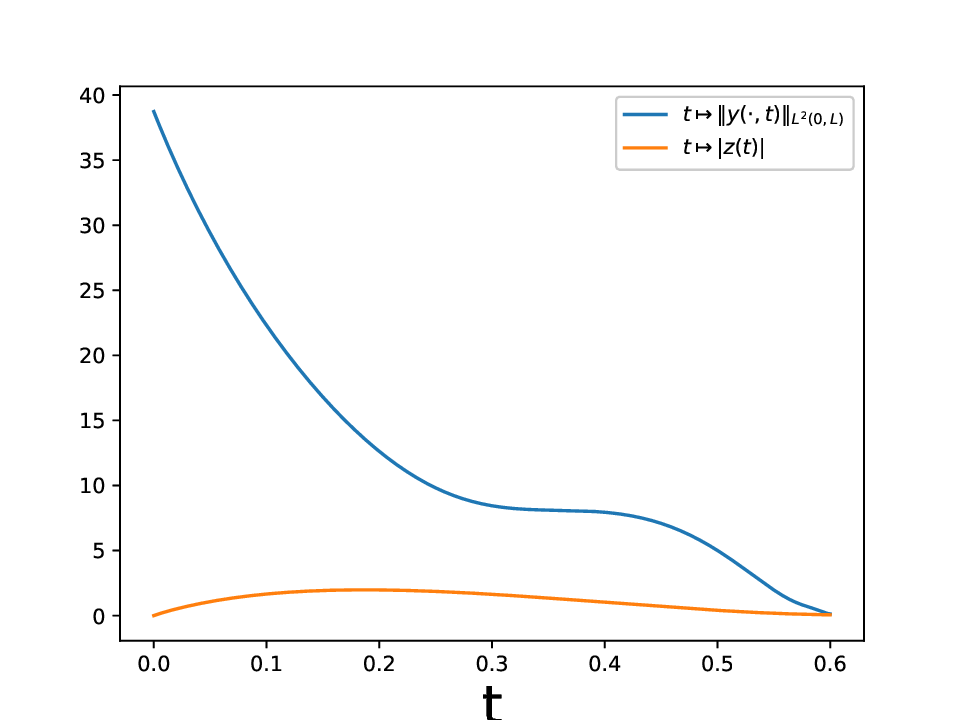}
		\end{center}
		\title{(a) Evolution of controlled state norms in time.}
	\end{minipage}
	\hfill
	\begin{minipage}[c]{.4\linewidth}
		\begin{center}
			\includegraphics[width=5cm]{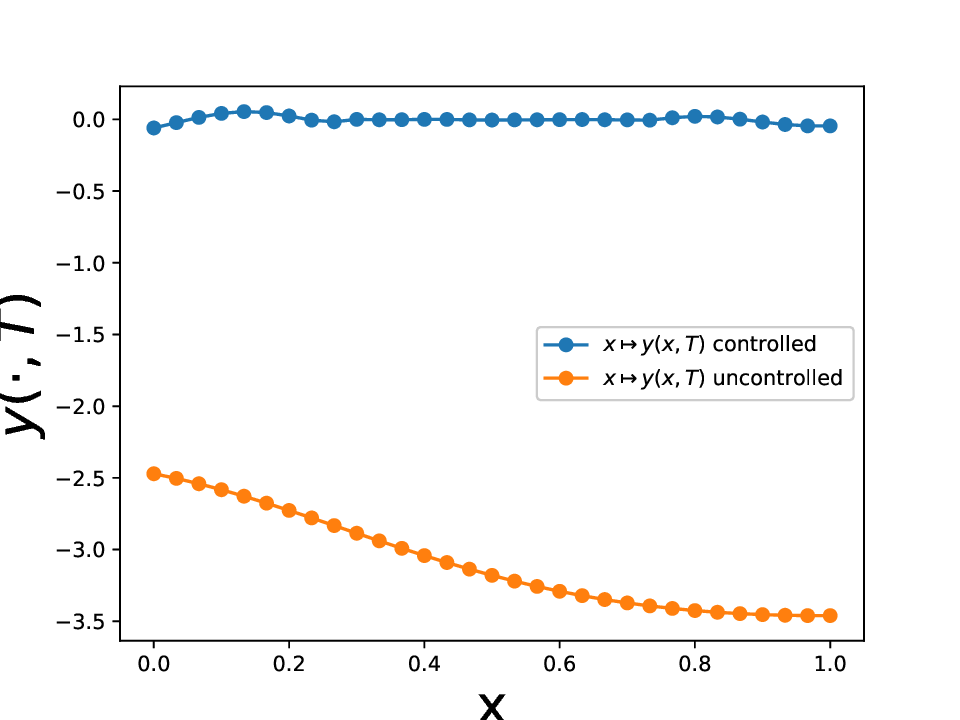}
		\end{center}
		\title{(b) Comparison of controlled and uncontrolled states at time $T$.}
	\end{minipage} 
	\caption{Evolution of controlled state norms in time and comparison of controlled and uncontrolled states at time $T$.}
\end{figure}

\begin{table}[h]
	\caption{Numerical results for Test 1.}\label{tab1}%
	\begin{tabular}{@{}cllll@{}}
		\toprule
		$\varepsilon$ &$10^{-1}$ & $10^{-2}$ & $10^{-3}$ & $10^{-4}$\\
		\midrule
		$N_{iter}$  &5 & 10 & 26 & 92 \\
		$\|y(\cdot,T)\|_{L^{2}(0,1)}$  & 1.4146 & 0.4614  & 0.099  & 0.0236 \\
		$|z(T)|$  & 1.8015 & 0.9726 & 0.2974  & 0.0598 \\
		$\|v\|_{L^{2}(\omega_T)}$  & 4.8074 & 9.2515  & 13.0852  & 14.9735  \\
		\botrule
	\end{tabular}
\end{table}

\subsubsection{Test 2} 
In a second experiment, we kept the data from Test 1, with the exception of the following. 
\begin{itemize}
	\item $b=1, y_0(x)=10\exp(-0.5(x-0.5)^2)$ and $z_0=y_0(0)$.
\end{itemize}
\begin{figure}[H]
	\begin{minipage}[c]{.3\linewidth}
		\begin{center}
			\includegraphics[width=4cm]{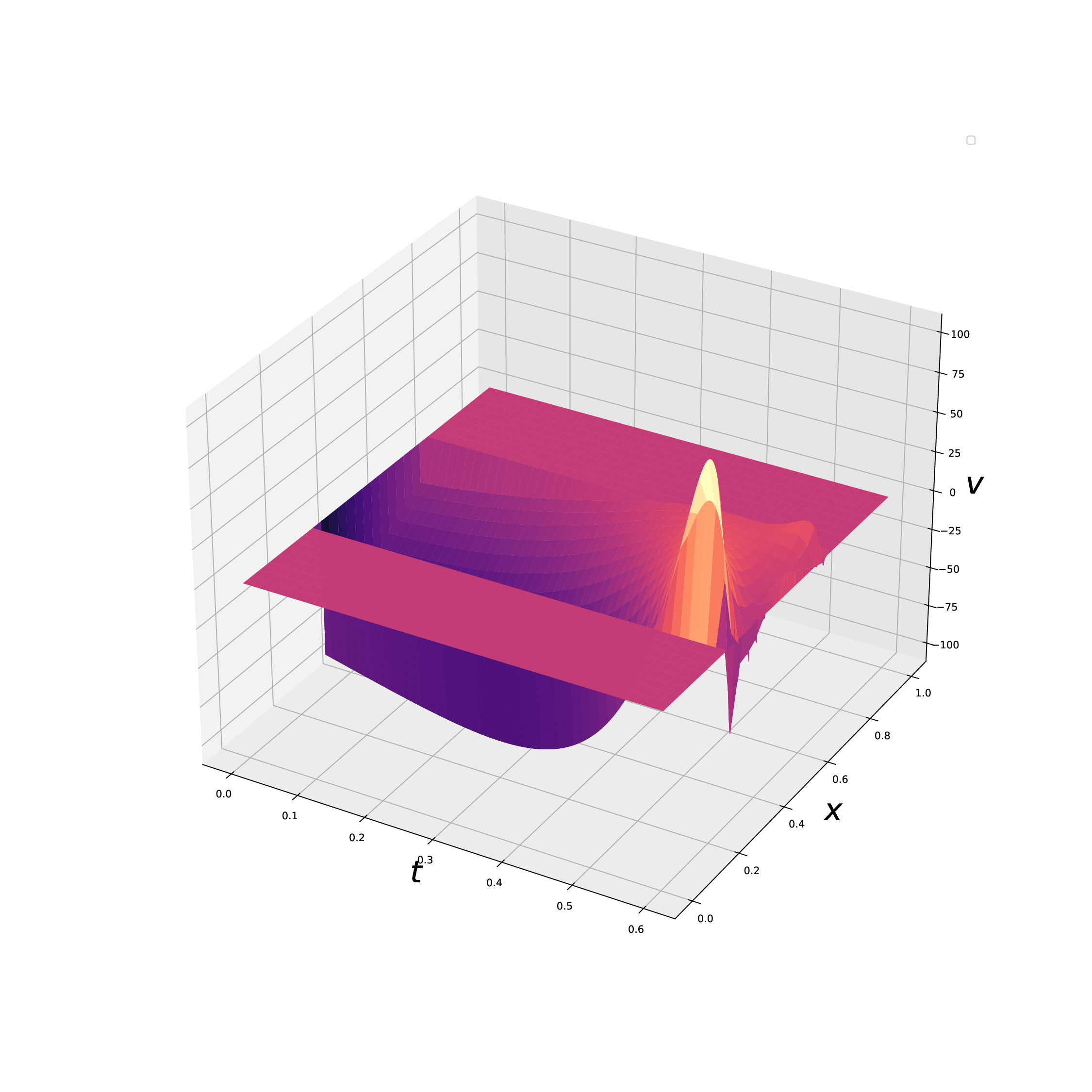}
		\end{center}
		\title{(a) Computed control.}
	\end{minipage} \hfill
	\begin{minipage}[c]{.3\linewidth}
		\begin{center}
			\includegraphics[width=4cm]{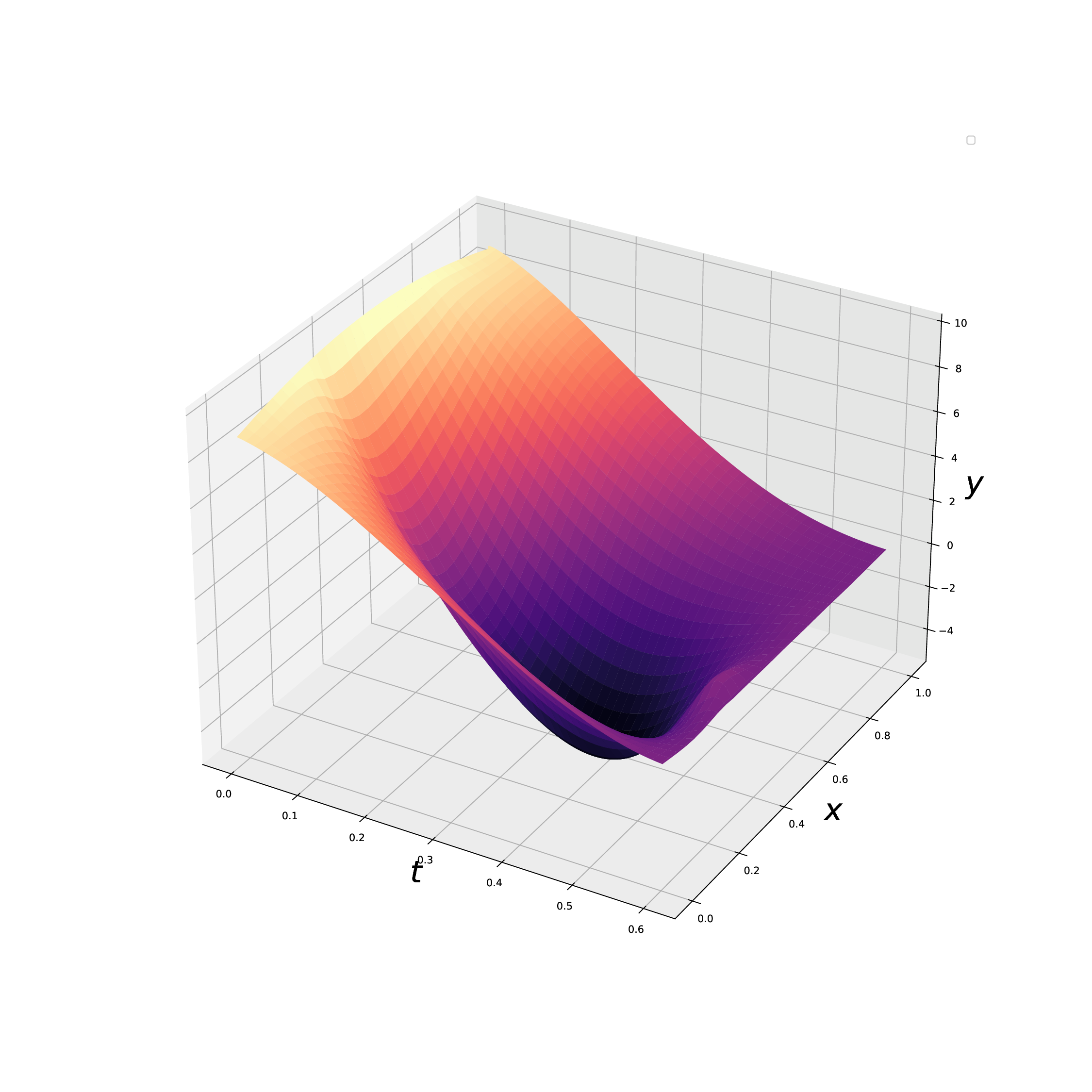}
		\end{center}
		\title{(b) Associated state $y$.}
	\end{minipage}
	\hfill
	\begin{minipage}[c]{.3\linewidth}
		\begin{center}
			\includegraphics[width=4cm]{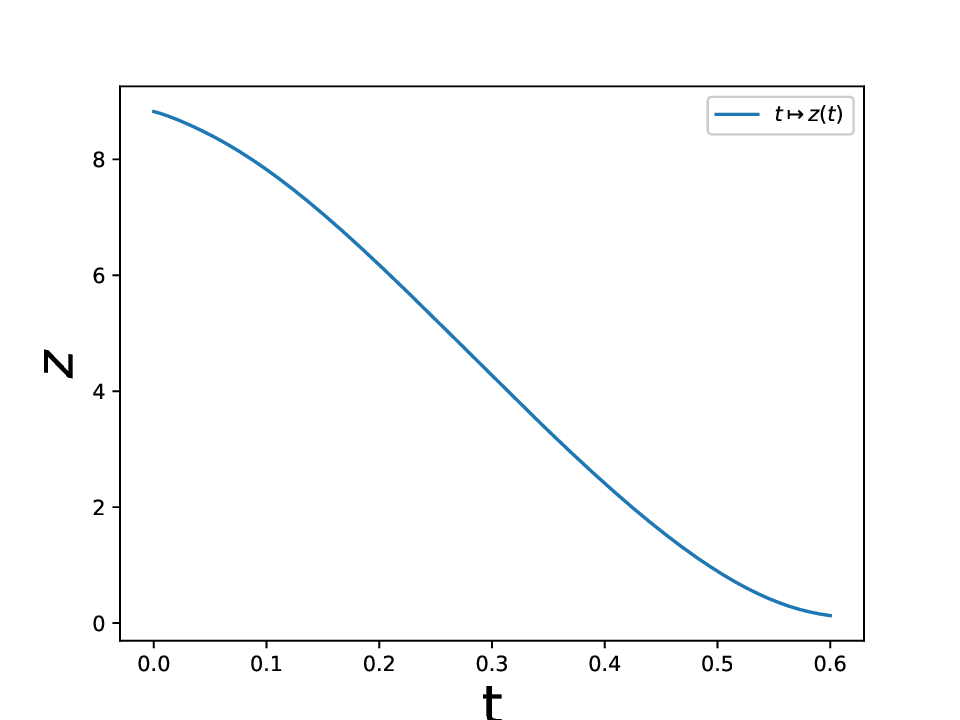}
		\end{center}
		\title{(c) Associated state $z$.}
	\end{minipage}
	\caption{The computed control and the associated state.}
\end{figure}

\begin{figure}[H]
	\begin{minipage}[c]{.5\linewidth}
		\begin{center}
			\includegraphics[width=5cm]{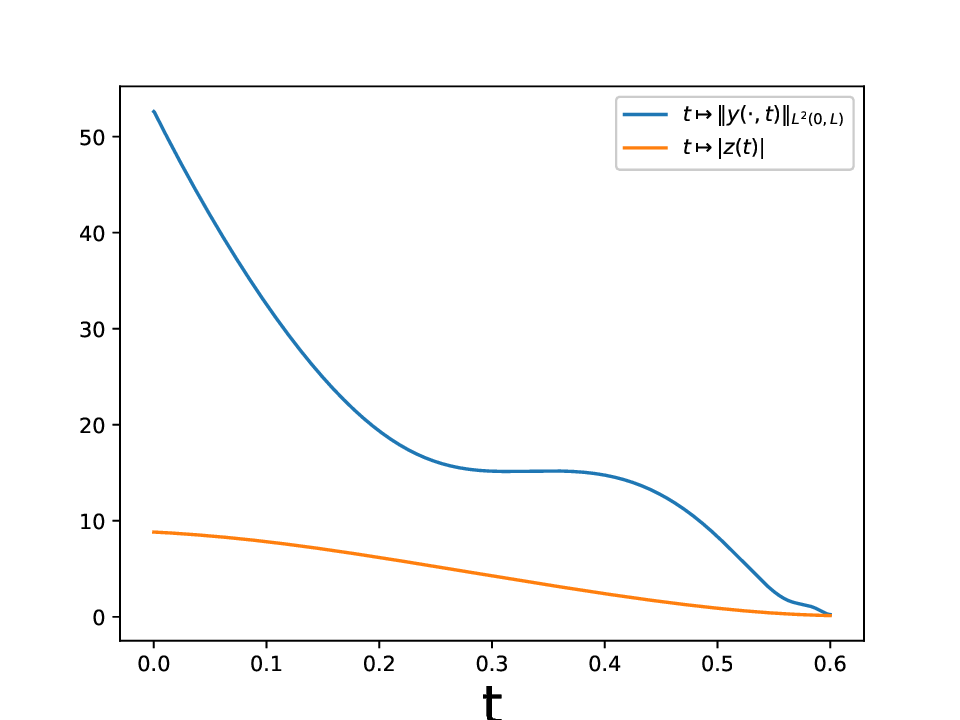}
		\end{center}
		\title{(a) Evolution of controlled state norms in time.}
	\end{minipage}
	\hfill
	\begin{minipage}[c]{.4\linewidth}
		\begin{center}
			\includegraphics[width=5cm]{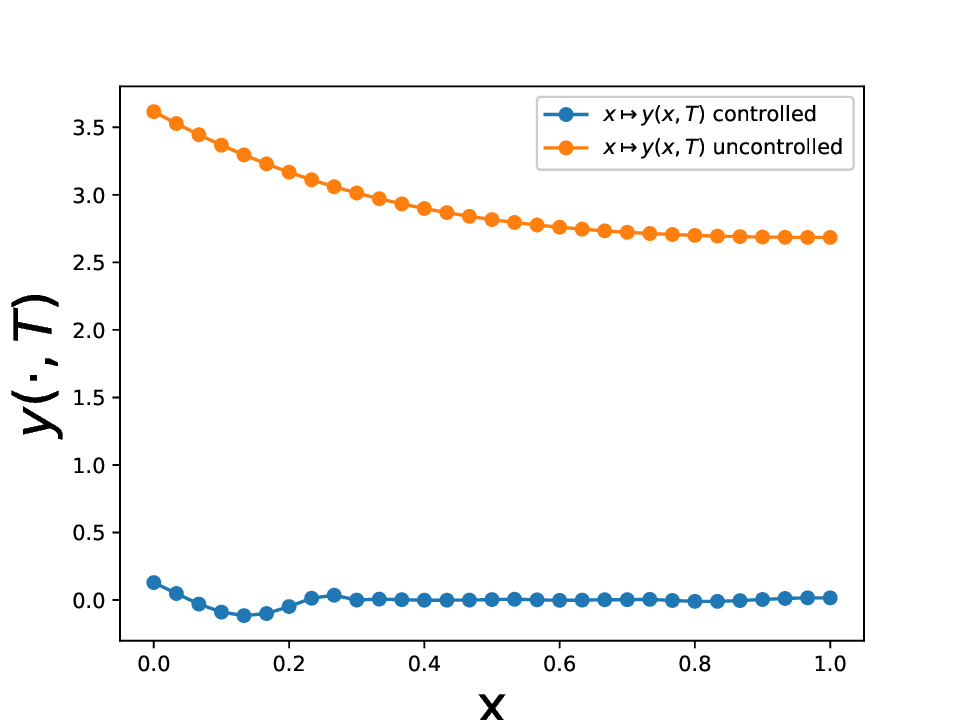}
		\end{center}
		\title{(b) Comparison of controlled and uncontrolled states at time $T$.}
	\end{minipage} 
	\caption{Evolution of controlled state norms in time and comparison of controlled and uncontrolled states at time $T$.}
\end{figure}

\begin{table}[h]
	\caption{Numerical results for Test 2.}\label{tab2}%
	\begin{tabular}{@{}cllll@{}}
		\toprule
		$\varepsilon$ &$10^{-1}$ & $10^{-2}$ & $10^{-3}$ & $10^{-4}$\\
		\midrule
		$N_{iter}$  &4 & 10 & 31 & 160 \\
		$\|y(\cdot,T)\|_{L^{2}(0,1)}$  & 1.5264 & 0.6296  & 0.1651  & 0.0431 \\
		$|z(T)|$  & 2.9587 & 1.8319 & 0.6150  & 0.1290 \\
		$\|v\|_{L^{2}(\omega_T)}$  & 4.7667 & 11.2303  & 19.5797  & 24.4130  \\
		\botrule
	\end{tabular}
\end{table}

\subsubsection{Test 3}
Finally, in a third experiment, we tried to analyze the problem of boundary control with the following data.
\begin{itemize}
	\item $\ell=1$, $T=0.6$.
	\item $y_{0}(x)=-10\sin(\pi x)$, $z_{0}=0$.
	\item $a=0, b=20, c=0$ and $\mu=\kappa=1$.
\end{itemize}
\begin{figure}[H]
	\begin{minipage}[c]{.3\linewidth}
		\begin{center}
			\includegraphics[width=4cm]{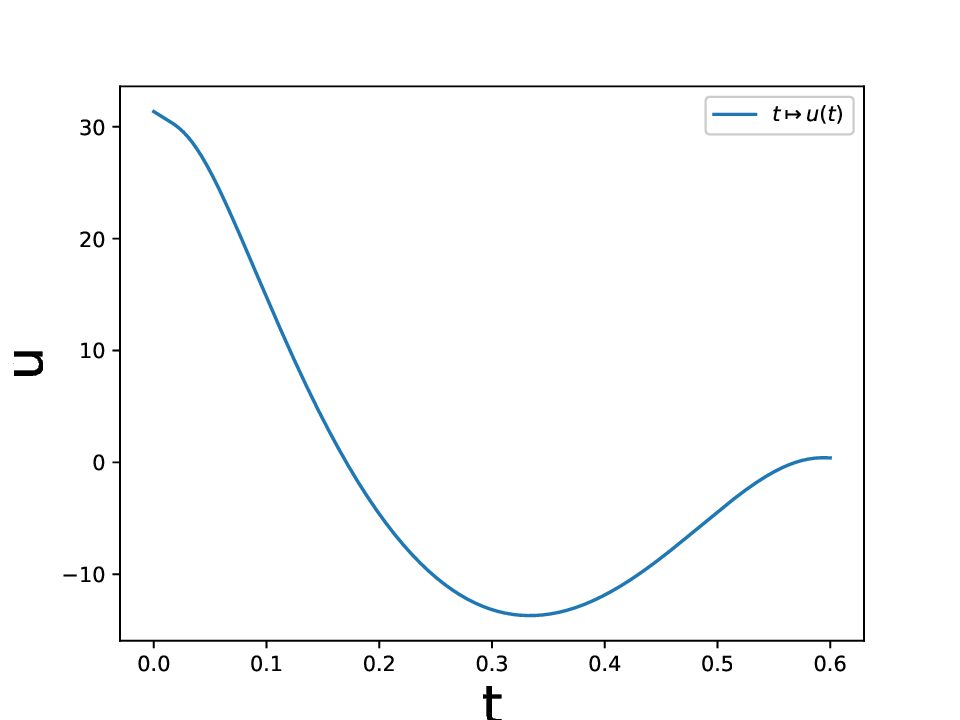}
		\end{center}.
		\title{(a) Computed control.}
	\end{minipage} \hfill
	\begin{minipage}[c]{.3\linewidth}
		\begin{center}
			\includegraphics[width=4cm]{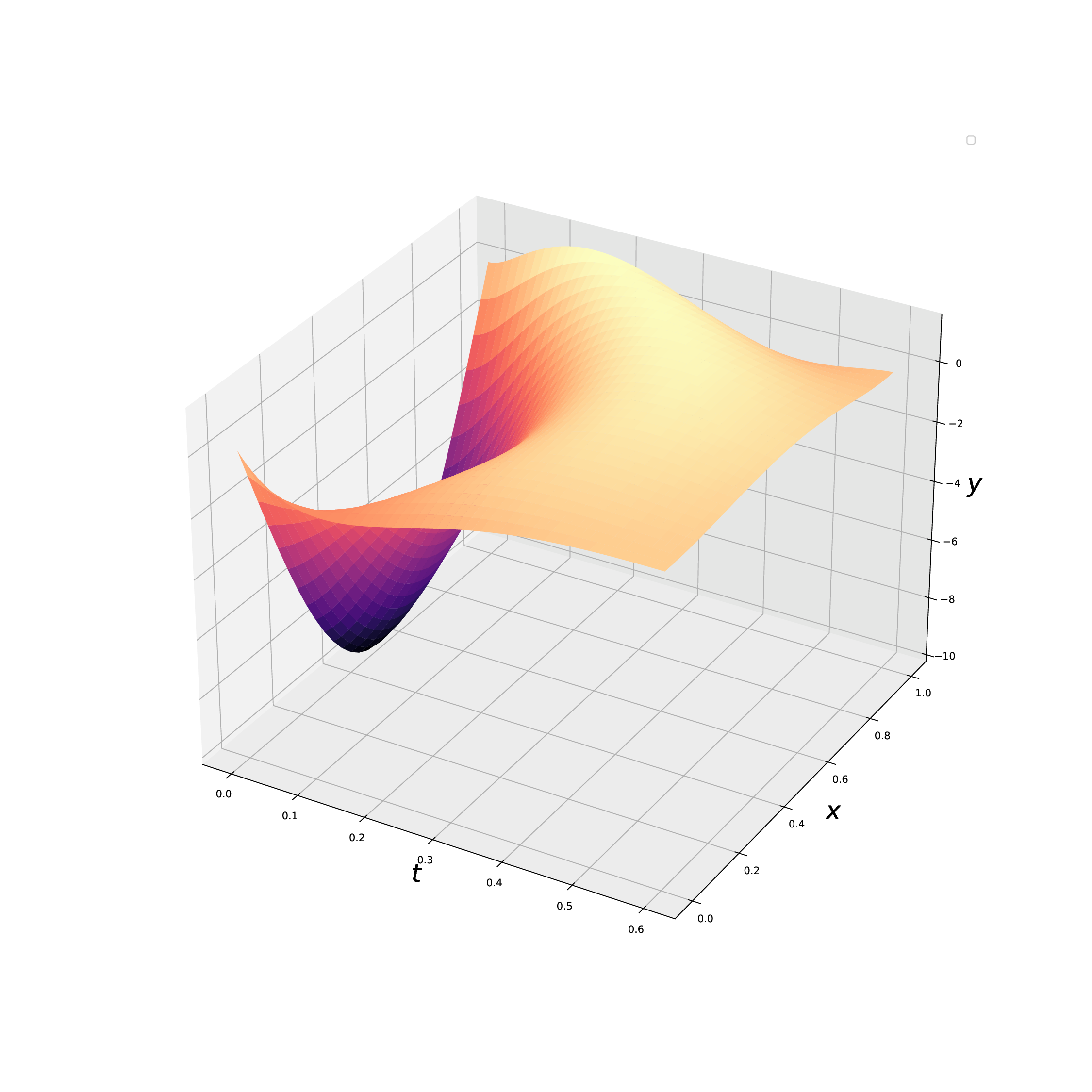}
		\end{center}
		\title{(b) Associated state $y$.}
	\end{minipage}
	\hfill
	\begin{minipage}[c]{.3\linewidth}
		\begin{center}
			\includegraphics[width=4cm]{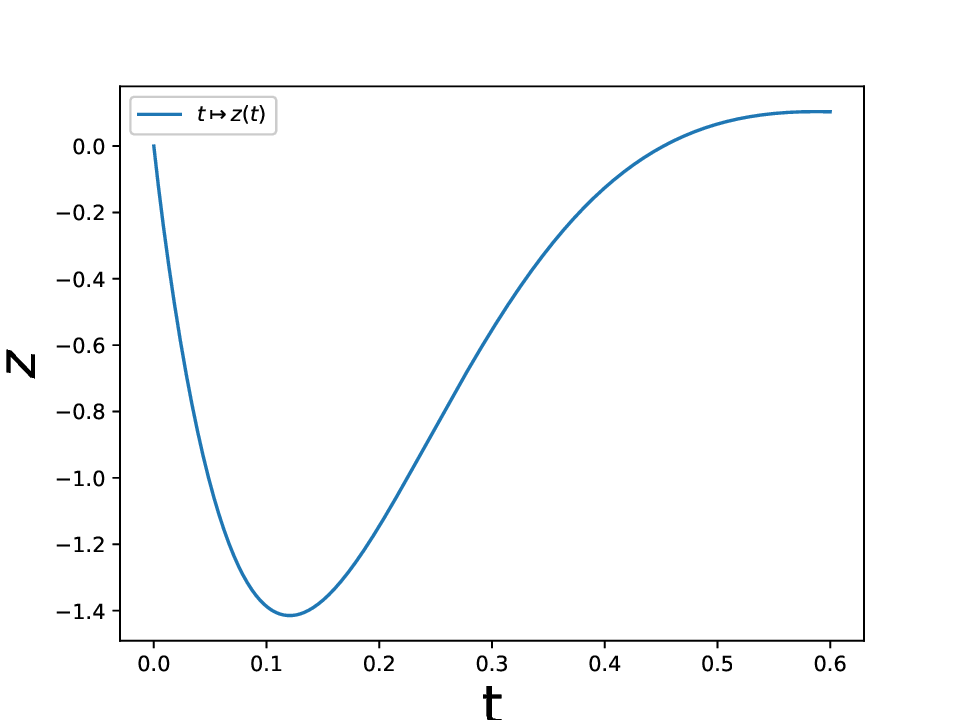}
			\title{(c) Associated state $z$.}
		\end{center}
	\end{minipage}
	\caption{The computed control and the associated state.}
\end{figure}

\begin{figure}[H]
	\begin{minipage}[c]{.5\linewidth}
		\begin{center}
			\includegraphics[width=5cm]{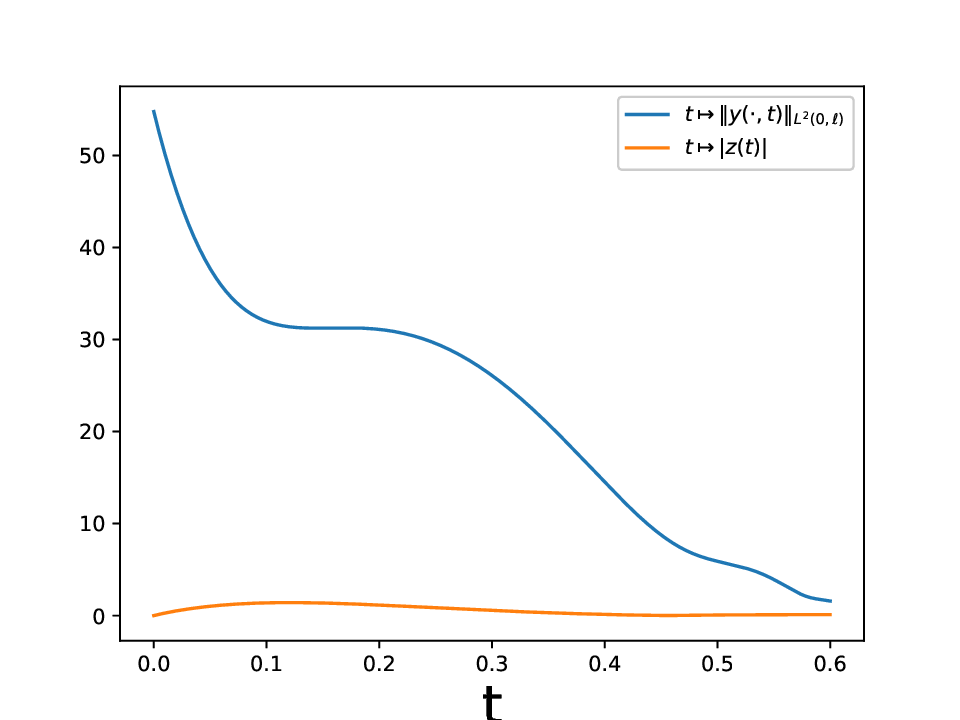}
		\end{center}
		\title{(a) Evolution of controlled state norms in time.}
	\end{minipage}
	\hfill
	\begin{minipage}[c]{.4\linewidth}
		\begin{center}
			\includegraphics[width=5cm]{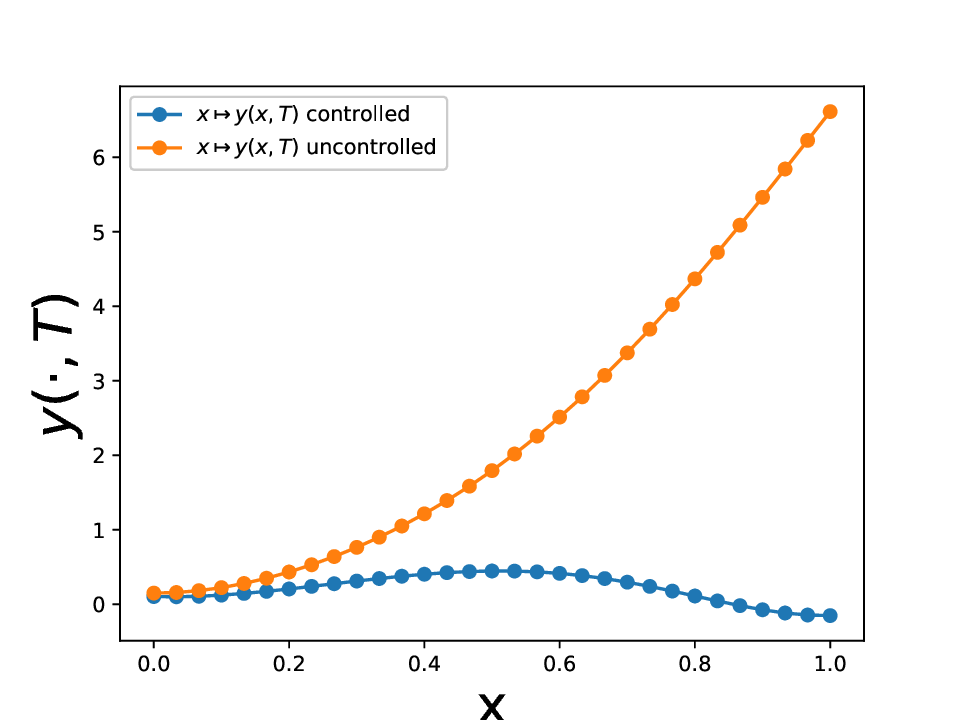}
			\title{(b) Comparison of controlled and uncontrolled states at time $T$.}
		\end{center}
	\end{minipage} 
	\caption{Evolution of controlled state norms in time and comparison of controlled and uncontrolled states at time $T$.}
\end{figure}

\begin{table}[h]
	\caption{Numerical results for Test 3.}\label{tab3}%
	\begin{tabular}{@{}cllll@{}}
		\toprule
		$\varepsilon$ &$10^{-1}$ & $10^{-2}$ & $10^{-3}$ & $10^{-4}$\\
		\midrule
		$N_{iter}$  &5 & 8 & 17 & 56 \\
		$\|y(\cdot,T)\|_{L^{2}(0,\ell)}$  & 1.9039 & 1.1385  & 0.5186  & 0.0431 \\
		$|z(T)|$  & 0.0410 & 0.0357 & 0.0283  & 0.0196 \\
		$\|v\|_{L^{2}(0,T)}$  & 12.4479 & 11.0431  & 10.3262  & 9.9543  \\
		\botrule
	\end{tabular}
\end{table}

In summary, numerical simulations of tests 1, 2 and 3 show that the HUM algorithm provides accurate results
for the numerical approximation of distributed or boundary controls for the heat equation coupled with an ordinary differential equation.
%\begin{table}[H]
%	\caption{Numerical results for Test 3.}\label{tab3}
%	\begin{tabular*}{\textwidth}{@{\extracolsep\fill}ccccc}
%		\toprule%
%		$\varepsilon$ &$10^{-1}$ & $10^{-2}$ & $10^{-3}$ & $10^{-4}$    \\
%		\midrule
%		$N_{iter}$  &5 & 8 & 17 & 56  \\
%		\midrule
%		$\|y(\cdot,T)\|_{L^{2}(0,\ell)}$  & 1.9039 & 1.1385  & 0.5186  & 0.0431  \\
%		\midrule
%		$|z(T)|$  & 0.0410 & 0.0357 & 0.0283  & 0.0196  \\
%		\midrule
%		$\|v\|_{L^{2}(0,T)}$  & 12.4479 & 11.0431  & 10.3262  & 9.9543  \\
%		\bottomrule
%	\end{tabular*}
%\end{table}

% Début de l'annexe
%\appendix
%\counterwithin{section}{appendix}

\section{Conclusion and future works}\label{sec13}
This paper investigates the null controllability of a coupled ODE-heat system with boundary control using the spatial domain extension method and Carleman estimates. First, the original system was transformed into an intermediate distributed control system with the same property of null controllability using a spatial domain extension. Next, we showed that the transformed system is well-posed and proved a new Carleman estimate for the corresponding adjoint system. Finally, the null controllability of the original control system was proved using the strategy of domain extension.
There are many interesting and important problems in this topic:
\begin{itemize}
	\item It is interesting to extend these results
	in higher dimensions: $y(x,t)\in\mathbb{R}^{n}, z(t)\in \mathbb{R}^{N}, a(x,t)\in\mathbb{R}^{N\times N}, b(x,t)\in\mathbb{R}^{n\times N}, c(t)\in\mathbb{R}^{N\times N}, \kappa\in \mathbb{R}^{N\times n}$ and $\mu \in \mathbb{R}^{n\times N}$. In our Carleman proof, we have simplified some terms thanks to the scalar product chosen as a function of $\mu$ and $\kappa$. In this case, such a choice is not obvious (at least for us). In the one-dimensional case $n=1$, $N\geq 1$, the authors of \cite{fernandez2010boundary} show null controllability results of \eqref{s0} for particular coefficients which do not depend on time thanks to the backstepping approach. But passing directly through an observability inequality for the adjoint system of \eqref{s0} remains an open question.
	\item We have proved that the ODE-heat system \eqref{s0} is null controllable with Neuman control and Dirichlet boundary coupling.  In the case of Dirichlet control and Dirichlet boundary coupling, we also achieve the same results with control in $H^{3/4}(0,T)$. It is also interesting to study the controllability properties with Dirichlet or Neumann control when the boundary coupling is of the Robin (or Fourier) type:
	\begin{equation*}  
		\left\{
		\begin{aligned}
			&y_{t}-y_{xx} + a(x,t)y+b(x,t)z(t)=0
			& & \text {in}\; Q_{\ell}, \\
			& z^{\prime}(t)+c(t)z(t)-\kappa y_{x}(0,t)
			=0 & & \text {in}\;(0,T), \\
			& y_{x}(0,t)+y(0,t)=\mu z(t)  & & \text {in}\;(0,T), \\
			& y(\ell,t) =u(t) & & \text {in}\;(0,T), \\
			& y(\cdot,0)=y_{0}  & & \text {in}\;(0,\ell), \\
			&z(0)=z_{0}.   \\	
		\end{aligned}
		\right.
	\end{equation*} 
\end{itemize}
\backmatter

\begin{appendices}
\section{Proof of Lemmas \ref{Lemma 1} and \ref{Lemma 2}} \label{Appendix A}
This paragraph is devoted to the proof of Lemmas \ref{Lemma 1} and \ref{Lemma 2}.\\
\begin{proof}[Proof of Lemma~{\upshape\ref{Lemma 1}}]
	Let $t\in (0,T)$ and $(y,\alpha), (\varphi,\beta)\in\mathcal{H}$.
	\begin{eqnarray*}
		\left\|\mathcal{B}(t)(y,\alpha) \right\|^{2}_{\mathcal{H}}&=&\int_{0}^{L}|a(x,t)y(x)+b(x,t)\alpha|^{2}dx+\frac{\mu}{\kappa}|c(t)|^{2}\alpha^{2}\\
		&\leq & 2\|a\|^{2}_{L^{\infty}(Q)}\|y\|^{2}_{L^{2}(0,L)}+2L\|b\|^{2}_{L^{2}(Q)}\alpha^{2}+\frac{\mu}{\kappa}\|c\|^{2}_{L^{\infty}(0,T)}\alpha^{2}\\
		&\leq & \left(2\|a\|^{2}_{L^{\infty}(Q)}+\frac{2L\kappa}{\mu}\|b\|^{2}_{L^{2}(Q)}+ \|c\|^{2}_{L^{\infty}(0,T)}\right)\left\|(y,\alpha) \right\|^{2}_{\mathcal{H}}.
	\end{eqnarray*}
	Consequently $\mathcal{B}(t)$ is uniformly bounded and we have $$\|\mathcal{B}(t)\|^{2}_{\mathcal{L}(\mathcal{H})}\leq 2\|a\|^{2}_{L^{\infty}(Q)}+2L\kappa\mu^{-1}\|b\|^{2}_{L^{2}(Q)}+ \|c\|^{2}_{L^{\infty}(0,T)}.$$ Now, we compute the adjoint of $\mathcal{B}(t)$.
	\begin{eqnarray*}
		\left\langle \mathcal{B}(t)(y,\alpha), (\varphi,\beta)\right\rangle_{\mathcal{H}}&=&\left\langle (
		-a(\cdot,t)y- b(\cdot,t)\alpha, - c(t)\alpha), (
		\varphi, \beta)\right\rangle_{\mathcal{H}}\\
		&=& -\int_{0}^{L}a(x,t)y(x)\varphi(x)dx -\alpha\int_{0}^{L}b(x,t)\varphi(x)dx-\frac{\mu}{\kappa} c(t)\alpha\beta \\
		&=& \left\langle (y,\alpha), \left(
		-a(\cdot,t)\varphi, -\frac{\kappa}{\mu}\displaystyle\int_{0}^{L}b(x,t)\varphi(x)dx-c(t)\beta \right)\right\rangle_{\mathcal{H}}.
	\end{eqnarray*}
	This provides the requested form of the adjoint of $\mathcal{B}(t)$.
\end{proof}
\begin{proof}[Proof of Lemma~{\upshape\ref{Lemma 2}}] Firstly, it is obvious that the form $\mathfrak{a}$ is well-defined, symmetric and positive.\\
	We claim that $\mathcal{H}^{1}$ is dense in $\mathcal{H}$. For that, it suffices to show that the orthogonal of $\mathcal{H}^{1}$ in $\mathcal{H}$ is trivial. Let $(y,\alpha)\in\mathcal{H}$ such that 
	\begin{eqnarray}
		\left\langle (y,\alpha),(
		\varphi, \beta)\right\rangle_{\mathcal{H}}=\int_{0}^{L}y(x)\varphi(x)dx+ \frac{\mu}{\kappa}\alpha\beta=0\quad \forall (
		\varphi, \beta) \in\mathcal{H}^{1}. \label{d1}
	\end{eqnarray}
	Choosing $\beta=0$ in \eqref{d1}, we obtain $\displaystyle\int_{0}^{L}y(x)\varphi(x)dx=0$, for all $\varphi\in D(0,L)$. Then $y=0$. Subsequently, since the trace operator is onto from $H^{1}(0,L)$ to $\mathbb{R}^{2}$, \eqref{d1} yields
	$$k^{-1}\alpha\varphi(0)=0\quad \forall \varphi\in H^{1}(0,L).$$
	Consequently $\alpha=0$.\\
	We claim that $\mathfrak{a}$ is closed. Let $(y,\alpha)\in\mathcal{H}^{1}$. It is obvious that 
	\begin{eqnarray*}
		\left\|(y,\alpha)\right\|^{2}_{\mathfrak{a}}:=\mathfrak{a}\left(  (y,\alpha),(y,\alpha)\right)+ \left\|(y,\alpha)\right\|^{2}_{\mathcal{H}}=\|y\|^{2}_{H^{1}(0,L)} + \frac{\mu}{\kappa}|\alpha|^{2}.
	\end{eqnarray*}
	Then $(D(\mathfrak{a}), \|\cdot\|^{2}_{\mathfrak{a}})$ is complete which is equivalent to the requested result (see \cite[Theorem 1.11]{kato2013perturbation}). 
\end{proof}

\section{Proof of the duality relation \eqref{dr}} \label{Appendix B}
By density, we can therefore restrict ourselves to initial values $(y_{0},z_{0})$ and final values $(\varphi_{T},\rho_{T})$ in $\mathcal{H}^{1}$, so that $(y,z)$ and $(
\varphi, \rho)$ are respictively strong solutions of systems \eqref{s1} and \eqref{s2}. Then
\begin{eqnarray*}
	\int_{\omega_{T}}v\varphi\d x\d t&=&  \int_{0}^{T}\langle (\mathds{1}_{\omega}v, 0), (\varphi(\cdot,t), \rho(t))\rangle_{\mathcal{H}}\d t \\
	&=& \int_{0}^{T}\langle \left[\partial_{t}+\mathcal{A}+\mathcal{B}(t)\right](y(\cdot,t),z(t)), (\varphi(\cdot,t), \rho(t))\rangle_{\mathcal{H}}\d t	\\
	&=& \int_{0}^{T}\langle (y(\cdot,t),z(t)), \underset{=0}{\underbrace{\left[-\partial_{t}+\mathcal{A}+(\mathcal{B}(t))^{*}\right](\varphi(\cdot,t), \rho(t))}}\rangle_{\mathcal{H}}\d t\\
	&& + \left[\langle (y(\cdot,t),z(t)), (\varphi(\cdot,t), \rho(t))\rangle_{\mathcal{H}}\right]_{t=0}^{t=T}\\
	&=& \left\langle (
	y(\cdot,T), z(T))
	, (\varphi_{T},\rho_{T}) \right\rangle_{\mathcal{H}}-\left\langle (
	y_{0}, z_{0}), (\varphi(\cdot,0), \rho(0))\right\rangle_{\mathcal{H}}.
\end{eqnarray*}
This shows the duality relation \eqref{dr}.
\end{appendices}
\bibliography{reffereces}% common bib file


\begin{thebibliography}{26}
\ifx \bisbn   \undefined \def \bisbn  #1{ISBN #1}\fi
\ifx \binits  \undefined \def \binits#1{#1}\fi
\ifx \bauthor  \undefined \def \bauthor#1{#1}\fi
\ifx \batitle  \undefined \def \batitle#1{#1}\fi
\ifx \bjtitle  \undefined \def \bjtitle#1{#1}\fi
\ifx \bvolume  \undefined \def \bvolume#1{\textbf{#1}}\fi
\ifx \byear  \undefined \def \byear#1{#1}\fi
\ifx \bissue  \undefined \def \bissue#1{#1}\fi
\ifx \bfpage  \undefined \def \bfpage#1{#1}\fi
\ifx \blpage  \undefined \def \blpage #1{#1}\fi
\ifx \burl  \undefined \def \burl#1{\textsf{#1}}\fi
\ifx \doiurl  \undefined \def \doiurl#1{\url{https://doi.org/#1}}\fi
\ifx \betal  \undefined \def \betal{\textit{et al.}}\fi
\ifx \binstitute  \undefined \def \binstitute#1{#1}\fi
\ifx \binstitutionaled  \undefined \def \binstitutionaled#1{#1}\fi
\ifx \bctitle  \undefined \def \bctitle#1{#1}\fi
\ifx \beditor  \undefined \def \beditor#1{#1}\fi
\ifx \bpublisher  \undefined \def \bpublisher#1{#1}\fi
\ifx \bbtitle  \undefined \def \bbtitle#1{#1}\fi
\ifx \bedition  \undefined \def \bedition#1{#1}\fi
\ifx \bseriesno  \undefined \def \bseriesno#1{#1}\fi
\ifx \blocation  \undefined \def \blocation#1{#1}\fi
\ifx \bsertitle  \undefined \def \bsertitle#1{#1}\fi
\ifx \bsnm \undefined \def \bsnm#1{#1}\fi
\ifx \bsuffix \undefined \def \bsuffix#1{#1}\fi
\ifx \bparticle \undefined \def \bparticle#1{#1}\fi
\ifx \barticle \undefined \def \barticle#1{#1}\fi
\bibcommenthead
\ifx \bconfdate \undefined \def \bconfdate #1{#1}\fi
\ifx \botherref \undefined \def \botherref #1{#1}\fi
\ifx \url \undefined \def \url#1{\textsf{#1}}\fi
\ifx \bchapter \undefined \def \bchapter#1{#1}\fi
\ifx \bbook \undefined \def \bbook#1{#1}\fi
\ifx \bcomment \undefined \def \bcomment#1{#1}\fi
\ifx \oauthor \undefined \def \oauthor#1{#1}\fi
\ifx \citeauthoryear \undefined \def \citeauthoryear#1{#1}\fi
\ifx \endbibitem  \undefined \def \endbibitem {}\fi
\ifx \bconflocation  \undefined \def \bconflocation#1{#1}\fi
\ifx \arxivurl  \undefined \def \arxivurl#1{\textsf{#1}}\fi
\csname PreBibitemsHook\endcsname

%%% 1
\bibitem[\protect\citeauthoryear{Krstic}{2009}]{krstic2009delay}
\begin{botherref}
\oauthor{\bsnm{Krstic}, \binits{M.}}:
Delay compensation for nonlinear, adaptive, and pde systems,
978--0
(2009)
\end{botherref}
\endbibitem

%%% 2
\bibitem[\protect\citeauthoryear{Fattorini and
  Russell}{1971}]{fattorini1971exact}
\begin{barticle}
\bauthor{\bsnm{Fattorini}, \binits{H.O.}},
\bauthor{\bsnm{Russell}, \binits{D.L.}}:
\batitle{Exact controllability theorems for linear parabolic equations in one
  space dimension}.
\bjtitle{Archive for Rational Mechanics and Analysis}
\bvolume{43}(\bissue{4}),
\bfpage{272}--\blpage{292}
(\byear{1971})
\end{barticle}
\endbibitem

%%% 3
\bibitem[\protect\citeauthoryear{Tenenbaum and
  Tucsnak}{2007}]{tenenbaum2007new}
\begin{barticle}
\bauthor{\bsnm{Tenenbaum}, \binits{G.}},
\bauthor{\bsnm{Tucsnak}, \binits{M.}}:
\batitle{New blow-up rates for fast controls of schr{\"o}dinger and heat
  equations}.
\bjtitle{Journal of Differential Equations}
\bvolume{243}(\bissue{1}),
\bfpage{70}--\blpage{100}
(\byear{2007})
\end{barticle}
\endbibitem

%%% 4
\bibitem[\protect\citeauthoryear{Lebeau and
  Robbiano}{1995}]{lebeau1995controle}
\begin{barticle}
\bauthor{\bsnm{Lebeau}, \binits{G.}},
\bauthor{\bsnm{Robbiano}, \binits{L.}}:
\batitle{Contr{\^o}le exact de l{\'e}quation de la chaleur}.
\bjtitle{Communications in Partial Differential Equations}
\bvolume{20}(\bissue{1-2}),
\bfpage{335}--\blpage{356}
(\byear{1995})
\end{barticle}
\endbibitem

%%% 5
\bibitem[\protect\citeauthoryear{Le~Rousseau and Lebeau}{2012}]{le2012carleman}
\begin{barticle}
\bauthor{\bsnm{Le~Rousseau}, \binits{J.}},
\bauthor{\bsnm{Lebeau}, \binits{G.}}:
\batitle{On carleman estimates for elliptic and parabolic operators.
  applications to unique continuation and control of parabolic equations}.
\bjtitle{ESAIM: Control, Optimisation and Calculus of Variations}
\bvolume{18}(\bissue{3}),
\bfpage{712}--\blpage{747}
(\byear{2012})
\end{barticle}
\endbibitem

%%% 6
\bibitem[\protect\citeauthoryear{Fursikov and
  Imanuvilov}{1996}]{fursikov1996controllability}
\begin{bbook}
\bauthor{\bsnm{Fursikov}, \binits{A.V.}},
\bauthor{\bsnm{Imanuvilov}, \binits{O.Y.}}:
\bbtitle{Controllability of Evolution Equations}.
\bpublisher{Seoul National University},
\blocation{Seoul}
(\byear{1996})
\end{bbook}
\endbibitem

%%% 7
\bibitem[\protect\citeauthoryear{Beauchard and
  Pravda-Starov}{2018}]{beauchard2018null}
\begin{barticle}
\bauthor{\bsnm{Beauchard}, \binits{K.}},
\bauthor{\bsnm{Pravda-Starov}, \binits{K.}}:
\batitle{Null-controllability of hypoelliptic quadratic differential
  equations}.
\bjtitle{Journal de l’{\'E}cole polytechnique-Math{\'e}matiques}
\bvolume{5},
\bfpage{1}--\blpage{43}
(\byear{2018})
\end{barticle}
\endbibitem

%%% 8
\bibitem[\protect\citeauthoryear{Miller}{2010}]{miller2010direct}
\begin{barticle}
\bauthor{\bsnm{Miller}, \binits{L.}}:
\batitle{A direct lebeau-robbiano strategy for the observability of heat-like
  semigroups}.
\bjtitle{Discrete and Continuous Dynamical Systems-Series B}
\bvolume{14}(\bissue{4}),
\bfpage{1465}--\blpage{1485}
(\byear{2010})
\end{barticle}
\endbibitem

%%% 9
\bibitem[\protect\citeauthoryear{Coron and Nguyen}{2017}]{coron2017null}
\begin{barticle}
\bauthor{\bsnm{Coron}, \binits{J.-M.}},
\bauthor{\bsnm{Nguyen}, \binits{H.-M.}}:
\batitle{Null controllability and finite time stabilization for the heat
  equations with variable coefficients in space in one dimension via
  backstepping approach}.
\bjtitle{Archive for Rational Mechanics and Analysis}
\bvolume{225},
\bfpage{993}--\blpage{1023}
(\byear{2017})
\end{barticle}
\endbibitem

%%% 10
\bibitem[\protect\citeauthoryear{Fern{\'a}ndez-Cara
  et~al.}{2006}]{fernandez2006null}
\begin{barticle}
\bauthor{\bsnm{Fern{\'a}ndez-Cara}, \binits{E.}},
\bauthor{\bsnm{Gonz{\'a}lez-Burgos}, \binits{M.}},
\bauthor{\bsnm{Guerrero}, \binits{S.}},
\bauthor{\bsnm{Puel}, \binits{J.-P.}}:
\batitle{Null controllability of the heat equation with boundary fourier
  conditions: the linear case}.
\bjtitle{ESAIM: Control, Optimisation and Calculus of Variations}
\bvolume{12}(\bissue{3}),
\bfpage{442}--\blpage{465}
(\byear{2006})
\end{barticle}
\endbibitem

%%% 11
\bibitem[\protect\citeauthoryear{Fern{\'a}ndez-Cara and
  Guerrero}{2006}]{fernandez2006global}
\begin{barticle}
\bauthor{\bsnm{Fern{\'a}ndez-Cara}, \binits{E.}},
\bauthor{\bsnm{Guerrero}, \binits{S.}}:
\batitle{Global carleman inequalities for parabolic systems and applications to
  controllability}.
\bjtitle{SIAM journal on control and optimization}
\bvolume{45}(\bissue{4}),
\bfpage{1395}--\blpage{1446}
(\byear{2006})
\end{barticle}
\endbibitem

%%% 12
\bibitem[\protect\citeauthoryear{Khoutaibi and
  Maniar}{2020}]{khoutaibi2020null}
\begin{barticle}
\bauthor{\bsnm{Khoutaibi}, \binits{A.}},
\bauthor{\bsnm{Maniar}, \binits{L.}}:
\batitle{Null controllability for a heat equation with dynamic boundary
  conditions and drift terms}.
\bjtitle{Evolution Equations \& Control Theory}
\bvolume{9}(\bissue{2}),
\bfpage{535}
(\byear{2020})
\end{barticle}
\endbibitem

%%% 13
\bibitem[\protect\citeauthoryear{Maniar et~al.}{2017}]{maniar2017null}
\begin{barticle}
\bauthor{\bsnm{Maniar}, \binits{L.}},
\bauthor{\bsnm{Meyries}, \binits{M.}},
\bauthor{\bsnm{Schnaubelt}, \binits{R.}}:
\batitle{Null controllability for parabolic equations with dynamic boundary
  conditions}.
\bjtitle{Evolution Equations \& Control Theory}
\bvolume{6}(\bissue{3}),
\bfpage{381}
(\byear{2017})
\end{barticle}
\endbibitem

%%% 14
\bibitem[\protect\citeauthoryear{Berinde et~al.}{2023}]{berinde2023qualitative}
\begin{botherref}
\oauthor{\bsnm{Berinde}, \binits{V.}},
\oauthor{\bsnm{Miranville}, \binits{A.}},
\oauthor{\bsnm{Moro{\c{s}}anu}, \binits{C.}}:
A qualitative analysis of a second-order anisotropic phase-field transition
  system endowed with a general class of nonlinear dynamic boundary conditions.
Discrete \& Continuous Dynamical Systems-Series S
\textbf{16}(1)
(2023)
\end{botherref}
\endbibitem

%%% 15
\bibitem[\protect\citeauthoryear{Fern{\'a}ndez-Cara
  et~al.}{2010}]{fernandez2010boundary}
\begin{barticle}
\bauthor{\bsnm{Fern{\'a}ndez-Cara}, \binits{E.}},
\bauthor{\bsnm{Gonz{\'a}lez-Burgos}, \binits{M.}},
\bauthor{\bsnm{Teresa}, \binits{L.}}:
\batitle{Boundary controllability of parabolic coupled equations}.
\bjtitle{Journal of Functional Analysis}
\bvolume{259}(\bissue{7}),
\bfpage{1720}--\blpage{1758}
(\byear{2010})
\end{barticle}
\endbibitem

%%% 16
\bibitem[\protect\citeauthoryear{Zeng et~al.}{2024}]{zeng2024null}
\begin{barticle}
\bauthor{\bsnm{Zeng}, \binits{C.}},
\bauthor{\bsnm{Zhou}, \binits{Z.}},
\bauthor{\bsnm{Xie}, \binits{C.}}:
\batitle{Null controllability of an ode-heat system coupled at boundary and
  internal term}.
\bjtitle{Applied Mathematics and Computation}
\bvolume{475},
\bfpage{128724}
(\byear{2024})
\end{barticle}
\endbibitem

%%% 17
\bibitem[\protect\citeauthoryear{Ammar-Khodja et~al.}{2011}]{ammar2011recent}
\begin{botherref}
\oauthor{\bsnm{Ammar-Khodja}, \binits{F.}},
\oauthor{\bsnm{Benabdallah}, \binits{A.}},
\oauthor{\bsnm{Gonz{\'a}lez~Burgos}, \binits{M.}},
\oauthor{\bsnm{Oteyza}, \binits{M.d.l.L.d.}}:
Recent results on the controllability of linear coupled parabolic problems: a
  survey.
Mathematical Control and Related Fields, 1 (3), 267-306.
(2011)
\end{botherref}
\endbibitem

%%% 18
\bibitem[\protect\citeauthoryear{Lions}{1988}]{lions1988controlabilite}
\begin{botherref}
\oauthor{\bsnm{Lions}, \binits{J.-L.}}:
Contr{\^o}labilit{\'e} exacte, stabilisation et perturbations de systemes
  distribu{\'e}s. tome 1. contr{\^o}labilit{\'e} exacte.
Rech. Math. Appl
\textbf{8}
(1988)
\end{botherref}
\endbibitem

%%% 19
\bibitem[\protect\citeauthoryear{Kato}{2013}]{kato2013perturbation}
\begin{bbook}
\bauthor{\bsnm{Kato}, \binits{T.}}:
\bbtitle{Perturbation Theory for Linear Operators}
vol. \bseriesno{132}.
\bpublisher{Springer},
\blocation{Berlin, Heidelberg}
(\byear{2013})
\end{bbook}
\endbibitem

%%% 20
\bibitem[\protect\citeauthoryear{Lunardi}{2009}]{Lunardi}
\begin{bbook}
\bauthor{\bsnm{Lunardi}, \binits{A.}}:
\bbtitle{Interpolation Theory}
vol. \bseriesno{Second edition}.
\bpublisher{Edizioni della Normale, Pisa},
\blocation{Pisa}
(\byear{2009})
\end{bbook}
\endbibitem

%%% 21
\bibitem[\protect\citeauthoryear{Lions}{2013}]{lions2013equations}
\begin{bbook}
\bauthor{\bsnm{Lions}, \binits{J.L.}}:
\bbtitle{Equations Differentielles Operationnelles: et Probl{\'e}mes aux
  Limites}
vol. \bseriesno{111}.
\bpublisher{Springer},
\blocation{Berlin, Heidelberg}
(\byear{2013})
\end{bbook}
\endbibitem

%%% 22
\bibitem[\protect\citeauthoryear{Pr{\"u}ss and
  Schnaubelt}{2001}]{pruss2001solvability}
\begin{barticle}
\bauthor{\bsnm{Pr{\"u}ss}, \binits{J.}},
\bauthor{\bsnm{Schnaubelt}, \binits{R.}}:
\batitle{Solvability and maximal regularity of parabolic evolution equations
  with coefficients continuous in time}.
\bjtitle{Journal of mathematical analysis and applications}
\bvolume{256}(\bissue{2}),
\bfpage{405}--\blpage{430}
(\byear{2001})
\end{barticle}
\endbibitem

%%% 23
\bibitem[\protect\citeauthoryear{Guerrero and
  Lebeau}{2007}]{guerrero2007singular}
\begin{barticle}
\bauthor{\bsnm{Guerrero}, \binits{S.}},
\bauthor{\bsnm{Lebeau}, \binits{G.}}:
\batitle{Singular optimal control for a transport-diffusion equation}.
\bjtitle{Communications in Partial Differential Equations}
\bvolume{32}(\bissue{12}),
\bfpage{1813}--\blpage{1836}
(\byear{2007})
\end{barticle}
\endbibitem

%%% 24
\bibitem[\protect\citeauthoryear{Lions and Magenes}{1972}]{lions1972non}
\begin{bbook}
\bauthor{\bsnm{Lions}, \binits{J.-L.}},
\bauthor{\bsnm{Magenes}, \binits{E.}}:
\bbtitle{Non-Homogeneous Boundary Value Problems and Applications, Volume II}.
\bpublisher{Springer},
\blocation{New York}
(\byear{1972})
\end{bbook}
\endbibitem

%%% 25
\bibitem[\protect\citeauthoryear{Glowinski et~al.}{2008}]{glowinskiexact}
\begin{botherref}
\oauthor{\bsnm{Glowinski}, \binits{R.}},
\oauthor{\bsnm{Lions}, \binits{J.-L.}},
\oauthor{\bsnm{He}, \binits{J.}}:
Exact and approximate controllability for distributed parameter systems: a
  numerical approach, 117.
Encyclopedia of mathematics and its applications
(2008)
\end{botherref}
\endbibitem

%%% 26
\bibitem[\protect\citeauthoryear{Boyer}{2013}]{boyer2013penalised}
\begin{bchapter}
\bauthor{\bsnm{Boyer}, \binits{F.}}:
\bctitle{On the penalised hum approach and its applications to the numerical
  approximation of null-controls for parabolic problems}.
In: \bbtitle{ESAIM: Proceedings},
vol. \bseriesno{41},
pp. \bfpage{15}--\blpage{58}
(\byear{2013}).
\bcomment{EDP Sciences}
\end{bchapter}
\endbibitem

\end{thebibliography}

%\section*{Declarations}
%$\bullet$ \textbf{Funding.}
%The authors declare that no funds, grants, or other support were received during the preparation of this manuscript.\\
%$\bullet$ \textbf{Competing Interests.}
%The authors have no relevant financial or non-financial interests to disclose.\\
%$\bullet$ \textbf{Author Contributions.} All authors contributed to the development of the paper, reviewed previous versions of the manuscript, and read and approved the final version.

\end{document}